\documentclass[final,3p,times,fleqn]{amsart}
\usepackage[margin=1in]{geometry}
\usepackage{amsaddr}

\input epsf
\usepackage{epsfig}
\usepackage{graphicx,xcolor}

\usepackage{fancyvrb}

\usepackage{url}
\usepackage[plainpages=false,pdfpagelabels,pdfencoding=auto,psdextra,colorlinks=true,urlcolor=black,linkcolor=black,citecolor=blue,hidelinks]{hyperref}
\usepackage{multirow}
\usepackage{booktabs}
\usepackage{epstopdf}
\usepackage{amsopn}
\usepackage{amsmath}
\usepackage{amsfonts}
\usepackage{amssymb}
\usepackage{graphicx}
\usepackage{epstopdf}
\usepackage{float}
\usepackage{mathtools}
\usepackage[inline]{enumitem}
\usepackage{bm}
\usepackage{underscore}
\usepackage[capitalize,nameinlink]{cleveref}
\usepackage{tikz}
\usetikzlibrary{patterns}

\usepackage{microtype}
\usepackage[normalem]{ulem}

\usepackage{caption}
\usepackage[labelformat=simple]{subcaption}

\usepackage{pgfplots}
\pgfplotsset{compat=1.17}
\usepackage{pgfplotstable}
\usetikzlibrary{matrix,backgrounds,decorations.pathreplacing,math}
\pgfdeclarelayer{myback}
\pgfsetlayers{myback,background,main}

\DeclareMathAlphabet\mathbold{OML}{cmm}{b}{it}
\numberwithin{equation}{section}

\newtheorem{remark}{Remark}[section]

\newtheorem{algorithm}{Algorithm}[section]

\renewcommand{\div}{\operatorname{div}}

\newcommand{\diag}{\operatorname{diag}}
\newcommand{\Var}{\operatorname{Var}}

\def\b1{{\mathbf 1}}

\newcommand{\R}{{\mathcal R}}

\newcommand{\C}{{\mathcal C}}

\newcommand{\T}{{\mathcal T}}

\renewcommand{\div}{\text{div}}

\newcommand{\hillary}[1]{\textcolor{black}{#1}}
\newcommand{\chak}[1]{{\textcolor{black}{#1}}}
\newcommand{\updates}[1]{\textcolor{black}{#1}}
\newcommand{\hillaryNEW}[1]{\textcolor{black}{#1}}

\title[]{Scalable Multilevel Monte Carlo Methods Exploiting Parallel Redistribution on Coarse Levels}
\author[]{Hillary R. Fairbanks$^{1}$ \and Delyan Z. Kalchev$^{1}$\and Chak Shing Lee$^{1}$\and Panayot S. Vassilevski$^{2}$}
\address{$^{1}$Center for Applied Scientific Computing, Lawrence Livermore National Laboratory, Livermore, CA, USA \\ $^{2}$Fariborz Maseeh Department of Mathematics and Statistics, Portland State University, Portland, OR, USA}

\date{\today}
\thanks{This work was performed under the auspices of the U.S. Department of Energy by Lawrence Livermore National Laboratory under contract DE-AC52-07NA27344 (LLNL-JRNL-836297).}

\begin{document}
\maketitle

\begin{abstract}
We study an element agglomeration coarsening strategy that requires data redistribution at coarse levels when the number of coarse elements becomes smaller than the used computational units (cores). The overall procedure generates coarse elements (general unstructured unions of fine grid elements) within the framework of element-based algebraic multigrid methods (or AMGe) studied previously. The AMGe generated coarse spaces have the ability to exhibit approximation properties of the same order as the fine-level ones since by construction they contain the piecewise polynomials  of the same order as the fine level ones. These approximation properties are key for the successful use of AMGe in multilevel solvers for nonlinear partial differential equations as well as for multilevel Monte Carlo (MLMC) simulations. The ability to coarsen without being constrained by the number of available cores, as described  in the present paper, allows to improve the scalability of these solvers as well as  in the overall MLMC method. The paper illustrates  {this latter fact with detailed scalability study of MLMC simulations} applied to model Darcy equations with a stochastic log-normal permeability field. 
\smallskip \\
\noindent \textbf{Keywords.} redistribution, distributed computing, multilevel methods, algebraic multigrid (AMG), AMGe, multilevel Monte Carlo
\end{abstract}


\section{Introduction}

In aggregation-based algebraic multigrid (AMG), when implemented in a distributed parallel environment, the coarsening process can reach a point when we need to aggregate current-level degrees of freedom (dofs) that are placed in different cores. 
For example, this can be necessary when communication costs between cores exceed the benefit of reducing the local problem size; but, this is especially the case when the number of new (coarse) aggregates becomes less than the number of cores.
In either case, we need to perform data redistribution to improve scaling on the coarser aggregated levels as well as allow for level-independent data distribution; that is, in multilevel approaches, the number of cores necessary on the finest level need not dictate the number of cores utilized on the coarser levels.
In this paper, we use an unstructured aggregation-based coarsening that enables us to perform
redistribution formulated in terms of parallel matrix operations (as described in full details in a more general setting in \cite{2022ParallelRedistributionReport}). We illustrate the performance of the approach in the setting of multilevel Monte Carlo (MLMC) simulations (\cite{Giles08, Heinrich01}) that exploits the element-based AMG method (commonly referred as AMGe) combined with the proposed redistribution. As it is well known for MLMC (as well as in multilevel methods for nonlinear partial differential equations (PDEs)), to achieve scalable performance we need to utilize a hierarchy of coarse spaces that exhibit guaranteed approximation properties in energy norm. These properties hold if we use hierarchy of spaces obtained by standard geometric refinement, however this is not necessarily the case for the more classical AMG methods (except in a modified two-level setting, see \cite{2019HuVassilevski}).   

In the application of MLMC, as well as in nonlinear multilevel solvers, for finite element or finite volume discretization problems posed in mixed saddle-point form, the ability of certain aggregation-based AMG to exhibit algorithmic scalability was demonstrated previously; see \cite{Barker21, Lee22}. The present paper extends these results to the case of element agglomeration AMG (AMGe) which utilizes general unstructured agglomerates incorporating, in addition, parallel redistribution coarsening. More specifically, the parallel scalability of the AMGe method is demonstrated within MLMC simulations for the Darcy equation with uncertain permeability field modeled by log-normal distribution (as commonly used in such studies, e.g., \cite{Cliffe11}, \cite{2021PosteriorMultilevel}).  We note that in the tests here,  we utilize the scalable PDE sampler utility from \cite{Osborn17, Osborn17b} (based on the work \cite{Lindgren11}) which contributes to the overall scalability of the MLMC method. 

We stress that the ability to coarsen beyond the number of cores is a key factor to achieve scalability of MLMC since additional coarser levels in the hierarchy results in fewer required fine solves -- a consequence of the variance reduction at finer levels. Therefore, in MLMC we benefit by having better scalable solvers at coarse levels. Incorporating redistribution to create as small as necessary coarse problems (and still possessing approximation quality as in AMGe) in addition to having the ability for scalable sampling is another key factor for achieving scalability of the overall MLMC process. \hillary{The former we view as the main contribution of the present paper, as it results in cost improvement for MLMC.}
We note that, while this work focuses on the parallelism of each PDE solve via
level-independent redistribution of mesh agglomerates, there are additional forms of parallelism from MLMC not included in this work. In particular, we do not employ scheduling algorithms to distribute (independent) simulations across the available cores, which is an ongoing area of research for large-scale MLMC (see, e.g.,~\cite{2017ParallelMLMC,shegunov2020dynamic,badia2021massively}).

The remainder of the paper is structured as follows. In Section~\ref{section: parallel finite elements}, we summarize some facts about finite element computations and introduce the main concept of parallel assembly utilizing  finite element  structures  such as element topology (elements, faces and edges) and associated degrees of freedom and element data. The involved algorithms  are formulated in terms of operations with parallel (distributed) sparse matrices. 
In Section~\ref{section: AMGe coarsening with redistribution}, we outline the particular AMGe coarsening (originated in  \cite{deRhamAMGe} and implemented as described in  \cite{2021HcurlHdivParELAG, 2021ParELAGreport}) utilizing element agglomeration, coarse element matrices and local and global interpolation matrices $P$ constructed elementwise and how they are modified after redistribution. 
Section~\ref{section: MLMC application} contains the application of the constructed AMGe hierarchy
of spaces on general unstructured agglomerates to MLMC. \hillaryNEW{We present detailed numerical results for the model Darcy equation of subsurface flow with uncertain permeability field in Section~\ref{section:numresults}, using a scalable PDE sampler technique (as in  \cite{Osborn17, Osborn17b}, and \cite{Barker21}). This section draws scaling and computational cost comparisons between a standard AMG hierarchy (without redistribution) with a hierarchy that incorporates our novel level-dependent redistribution. These results, which include scalability analysis and MLMC comparisons show that the level-dependent redistribution is able to improve scaling and total CPU cost (online and offline) without sacrificing accuracy.} Finally, in Section~\ref{section: conclusions}, we draw some conclusions and outline some directions for  possible future extensions.

\section{Parallel finite element assembly}\label{section: parallel finite elements}
Here, we summarize some facts about finite element computations in the form presented in \cite{2002AMGeTopology} (see also \cite{VassilevskiMG}) before extending them in the setting of element redistribution in the following section. 

We are given a finite element mesh $\T_h = \{\tau\}$ and a bilinear form $a(.,.)$ that is a sum of local element-based bilinear forms $a_\tau(.,.)$. The parallel data, which is initially distributed across all available cores, is represented through matrix-based \emph{relations}, such as ``element\_dof'' and
``dof\_truedof"\updates{, reviewed below}. Here the former relation is essential to map the element matrices to the corresponding dofs, and the latter assembles the global system from element matrices in parallel. These matrices -- ``object1\_object2'' -- \chak{are} boolean sparse matrices (i.e., matrices with entries equal to either 0 or 1) where object1's entries are listed as rows and object2's entries are listed as columns. For ``element\_dof", the matrix has nonzero entries at positions $(i,j)$ if element $i$ contains dof $j$. \chak{Figure~\ref{fig:relation-table-example-a}-\subref{fig:relation-table-example-b} shows an example of a simple 2D mesh with a set of elements and dofs of the Raviart-Thomas finite element space, while the corresponding element\_dof relation is given in Figure~\ref{fig:element-dof}.} Note, the transpose relation ``dof\_element" $=($``element\_dof"$)^T$ shows in its $j$th row all elements $i$ that dof $j$ belongs to. Furthermore, using matrix operations, one can create composite relations. For example, given the relation ``element\_\chak{facet}" representing each element in terms of its \chak{facets (a facet means a face in 3D meshes or an edge in 2D meshes)}, we can form its transpose
\[
\text{``\chak{facet}\_element"=$(``$element\_\chak{facet}$")^T$}
\]
 and by forming the product ``element\_element" = ``element\_\chak{facet}" $\times$``\chak{facet}\_element", we get the new relation that relates two elements if they share a common \chak{facet}. Such a relation is commonly used as an input graph for software producing mesh/graph partitioners. 

In a typical finite element setting, in addition to the element topology represented as relation tables like ``element\_\chak{facet}" and ``\chak{facet\_vertex}",
one needs the relation ``element\_dof", mentioned above and further discussed below. The latter depends on the finite element spaces used to evaluate the bilinear form
$a(.,.)$  which is typically done on an element-by-element manner utilizing the element matrices $a_\tau(.,.)$. 
In a parallel computing environment, for our purpose, it is convenient to consider the degrees of freedom that are decoupled in the following sense. 
\chak{In the case of $H(\div)$-conforming finite elements (lowest-order case), each true degree of freedom (referred to as ``truedof") is associated with one and only one facet of the mesh. 
We assign to the same facet a separate degree of freedom for each of the elements sharing the facet.
These (decoupled) degrees of freedom are referred to as ``dofs". Figure~\ref{fig:relation-table-example-a}-\subref{fig:relation-table-example-b} is an example demonstrating the connection between ``truedofs" and ``dofs".
In our implementation, the relation of each ``truedof" and ``dof" is stored in the boolean matrix ``dof\_truedof" (see Figure~\ref{fig:truedof-dof}).
}

\updates{In general, true dofs are defined as the uniquely identifiable global dofs coming from the mathematical definition of the respective finite element space and associated with global basis functions. That is, true dofs are not modified by particular implementation details such as cloning or decoupling due to local element considerations or parallel distribution and inter-core sharing of the dofs (sometimes referred as ghost dofs). In connection to that, a mapping relation like ``dof\_truedof" is helpful in linking decoupled, cloned, or ghost dofs back to their respective unique global true dofs. In the lowest-order $H(\div)$-conforming example (see Figure~\ref{fig:relation-table-example-a},\subref{fig:truedof-dof}), the true dofs are uniquely associated with global facets independently of the adjacent elements or if a facet is at the interface between two cores.}

In what follows, we assume that we are given the relation table ``element\_dof" which relates every element with the (decoupled) dofs. 
We can form also the composite relation ``dof\_dof" $=($``element\_dof"$)^T\times$ ``element\_dof" which couples all dofs within each element but not across the elements. 
Finally, we can compute the element matrices $A_\tau = \left (a_\tau(\phi_j,\;\phi_i) \right )_{i,j \text{ are dofs in }\tau}$, with basis functions $\phi_i, \phi_j$. This is typically done on  a reference element.  
Since each element enters in a relation ``element\_dof" where the dofs are listed in a particular global order, that order, e.g.,  $j_1,j_2,\dots, \;j_m$ also implies a specific local (to the element) order $1,2,\;\dots, m$\updates{, providing the basic order-preserving index mapping ($j_k \leftrightarrow k$) between a global dof $j_k$ and a local dof $k$}. For example, for triangular elements (lowest order) we have $m=3$. This means that $A_\tau$ is a small $m \times m$ dense matrix. 
We can form a global matrix with \updates{a sparsity pattern} from the relation 
``dof\_dof" by placing at position \updates{$(j_k,j_l)$} the respective entry \updates{$(A_\tau)_{k,l}$, utilizing the basic index mapping between local and global numbering}. The observation is that the resulting matrix $A_{\diag}$ is block-diagonal with blocks  corresponding to each element $\tau$,  and  diagonal blocks equal to $A_\tau$.  

The actual assembled  matrix $A$, that we are interested in solving problems with, is obtained by performing numerically  the product $A = P^T A_{\diag} P$ where $P$ is the boolean matrix 
\updates{representing} the relation ``dof\_truedof". The sparsity pattern of $A$ corresponds to the pattern of the composite relation 
\[
\text{``truedof\_truedof" $= ($``dof\_truedof"$)^T \times$ ``dof\_dof" $\times$ ``dof\_truedof".}
\]

We note that all operations above  can be performed in a distributed parallel computing environment as long as the relations and associated matrices are stored in a distributed sparse matrix format.

\section{On the AMGe coarsening combined with redistribution}\label{section: AMGe coarsening with redistribution} 
Next, we give some brief details about element agglomeration coarsening within AMGe. Assume that we have generated a relation ``[agglomerated element]\_element" or ``AE\_element" for short,
which is a partitioning of the set of elements $\{\tau\}$  into subdomains $\{T\}$ where each subdomain $T$ is a connected union of elements $\tau$ \chak{(see Figure~\ref{fig:relation-table-example} for an illustration)}. 
\chak{For example, a partitioning can be generated using} a mesh partitioner (like METIS, \cite{metis}) utilizing \chak{the} relation (graph) ``element\_element" \chak{described in Section~\ref{section: parallel finite elements}}.

Then the set of agglomerated elements (AEs)  is used as a set of coarse elements. By local AE-by-AE computations  we build local but conforming interpolation matrices $\{P_T\}$ for each AE $T$. 
Conforming refers to the fact that if two AEs, $T_1$ and $T_2$, share a common true dof $i$ \updates{on the fine level}, the column indices of \updates{the nonzero entries} of $P_{T_1}$ and $P_{T_2}$ in their \updates{respective rows mapped to the true dof $i$} correspond to the same coarse true dofs, and the corresponding entry values coincide. This allows to define a unique global $P$ coming from the local AE {operators} $\{P_T\}$. Additionally, one is able to define coarse element matrices $A^c_T = P^T_T A_T P_T$,
 where $A_T$ is assembled from $\{A_\tau\}_{\tau \subset T}$. Finally, the conformity ensures that the global coarse matrix $A^c = P^T A P$ is the same as the assembled one from the local
coarse element matrices $\{A^c_T\}$. For more details we refer to \cite{VassilevskiMG}, 
 \cite{deRhamAMGe}, and its parallel implementation in \cite{2021HcurlHdivParELAG, 2021ParELAGreport}.

\def\elem{{e}}
\def\dof{{d}}
\def\truedof{{td}}
\tikzmath{ \unit = 7; }

 \begin{figure}
	\small
	\centering
	\begin{subfigure}[b]{0.48\textwidth}
		\centering
		\begin{tikzpicture}
		[
		cellCentroid/.style={circle,draw=black,text=black,fill=black,thick,minimum size=1mm, inner sep = 0em, outer sep = 0.0em},
		edgeCentroid/.style={rectangle,text=black,thick,minimum size=1mm, inner sep = 0em, outer sep = 0.0em},
		cellLabel/.style={circle,draw=black,text=black,fill=red!50!black!20,thick,minimum size=4mm, inner sep = 0.1em, outer sep = 0.1em},
		edgeLabel/.style={rectangle,draw=black,text=black,fill=blue!50!black!20,thick,minimum size=4mm, inner sep = 0em, outer sep = 0.1em},
		AE1/.style={rectangle,draw=blue!50!black!20,fill=blue!50!black!20,thick,minimum size=4mm, inner sep = 0em, outer sep = 0.0em},
		AE2/.style={rectangle,draw=red!50!black!20,fill=red!50!black!20,thick,minimum size=4mm, inner sep = 0em, outer sep = 0.0em}]
		
		\tikzmath{ \x1 = 0.0*\unit; \x2 = 0.25*\unit; \x3 = 0.5*\unit; \x4 = 0.75*\unit; \x5 = 0.95*\unit; }

		\coordinate [] (v1) at ( \x1, \x1);
		\coordinate [] (v2) at ( \x2, \x1);
		\coordinate [] (v3) at ( \x3, \x1);
		\coordinate [] (v4) at ( \x4, \x1);
		\coordinate [] (v6) at ( \x1, \x2);
		\coordinate [] (v7) at ( \x2, \x2);
		\coordinate [] (v8) at ( \x3, \x2);
		\coordinate [] (v9) at ( \x4, \x2);
		\coordinate [] (v11) at ( \x1, \x3);
		\coordinate [] (v12) at ( \x2, \x3);
		\coordinate [] (v13) at ( \x3, \x3);
		\coordinate [] (v14) at ( \x4, \x3);
		\coordinate [] (v15) at ( \x5, \x2);
		\coordinate [] (v16) at ( \x5, \x3);

		\draw[draw=blue!50!black!20,fill=blue!50!black!20] (v1) -- (v2) -- (v7) -- (v6) -- cycle;
		\draw[draw=red!50!black!20,fill=red!50!black!20] (v2) -- (v3) -- (v8) -- (v7) -- cycle;
		\draw[draw=red!50!black!20,fill=red!50!black!20] (v3) -- (v4) -- (v9) -- (v8) -- cycle;
		\draw[draw=blue!50!black!20,fill=blue!50!black!20] (v6) -- (v7) -- (v12) -- (v11) -- cycle;
		\draw[draw=red!50!black!20,fill=red!50!black!20] (v7) -- (v8) -- (v13) -- (v12) -- cycle;
		\draw[draw=red!50!black!20,fill=red!50!black!20] (v8) -- (v9) -- (v14) -- (v13) -- cycle;

		\draw[] (v1) -- (v4) -- (v14) -- (v11) -- cycle;
		\draw[] (v6) -- (barycentric cs:v6=0.65,v7=0.35);
		\draw[] (barycentric cs:v6=0.35,v7=0.65) -- (barycentric cs:v7=0.65,v8=0.35);
		\draw[] (barycentric cs:v7=0.35,v8=0.65) -- (barycentric cs:v8=0.65,v9=0.35);
		\draw[] (barycentric cs:v8=0.35,v9=0.65) -- (v9);

		\draw[] (v2) -- (barycentric cs:v2=0.6,v7=0.4);
		\draw[] (barycentric cs:v2=0.37,v7=0.63) -- (barycentric cs:v7=0.6,v12=0.4);
		\draw[] (barycentric cs:v7=0.37,v12=0.63) -- (v12);

		\draw[] (v3) -- (barycentric cs:v3=0.6,v8=0.4);
		\draw[] (barycentric cs:v3=0.37,v8=0.63) -- (barycentric cs:v8=0.6,v13=0.4);
		\draw[] (barycentric cs:v8=0.37,v13=0.63) -- (v13);

		\node[] at (barycentric cs:v1=0.5,v7=0.5) [] {$\elem_1$};
		\node[] at (barycentric cs:v6=0.5,v12=0.5) [] {$\elem_2$};
		\node[] at (barycentric cs:v2=0.5,v8=0.5) [] {$\elem_3$};
		\node[] at (barycentric cs:v7=0.5,v13=0.5) [] {$\elem_4$};
		\node[] at (barycentric cs:v3=0.5,v9=0.5) [] {$\elem_5$};
		\node[] at (barycentric cs:v8=0.5,v14=0.5) [] {$\elem_6$};
		
		\node[label=west:$AE_1:$ ] at (barycentric cs:v15=0.3,v16=0.7) [AE1] {};
		\node[label=west:$AE_2:$ ] at (barycentric cs:v15=0.65,v16=0.35) [AE2] {};

		\node[] at (barycentric cs:v2=0.5,v7=0.5) [edgeCentroid] {\scriptsize$\truedof_1$};
		\node[] at (barycentric cs:v3=0.5,v8=0.5) [edgeCentroid] {\scriptsize$\truedof_4$};
		\node[] at (barycentric cs:v7=0.5,v12=0.5) [edgeCentroid] {\scriptsize$\truedof_2$};
		\node[] at (barycentric cs:v8=0.5,v13=0.5) [edgeCentroid] {\scriptsize$\truedof_5$};
		\node[] at (barycentric cs:v6=0.5,v7=0.5) [edgeCentroid] {\scriptsize$\truedof_3$};
		\node[] at (barycentric cs:v7=0.5,v8=0.5) [edgeCentroid] {\scriptsize$\truedof_6$};
		\node[] at (barycentric cs:v8=0.5,v9=0.5) [edgeCentroid] {\scriptsize$\truedof_7$};

		\end{tikzpicture}

		\caption{Enumeration of elements ($\elem_i$), true dofs ($\truedof_i$), and AEs.}
		\label{fig:relation-table-example-a}
	\end{subfigure}
	\hfill
	\begin{subfigure}[b]{0.48\textwidth}
		\centering
		\begin{tikzpicture}
		[
		cellCentroid/.style={circle,draw=black,text=black,fill=black,thick,minimum size=1mm, inner sep = 0em, outer sep = 0.0em},
		edgeCentroid/.style={rectangle,text=black,thick,minimum size=1mm, inner sep = 0em, outer sep = 0.0em},
		cellLabel/.style={circle,draw=black,text=black,fill=red!50!black!20,thick,minimum size=4mm, inner sep = 0.1em, outer sep = 0.1em},
		edgeLabel/.style={rectangle,draw=black,text=black,fill=blue!50!black!20,thick,minimum size=4mm, inner sep = 0em, outer sep = 0.1em},
		AE1/.style={rectangle,draw=blue!50!black!20,fill=blue!50!black!20,thick,minimum size=4mm, inner sep = 0em, outer sep = 0.0em},
		AE2/.style={rectangle,draw=red!50!black!20,fill=red!50!black!20,thick,minimum size=4mm, inner sep = 0em, outer sep = 0.0em}]
		
		\tikzmath{
			\x1 = 0.0*\unit; \x2 = 0.25*\unit; \x3 = 0.5*\unit; \x4 = 0.75*\unit; }

		\coordinate [] (v11) at ( \x1, \x1);
		\coordinate [] (v12) at ( \x2, \x1);
		\coordinate [] (v13) at ( \x3, \x1);
		\coordinate [] (v14) at ( \x4, \x1);
		\coordinate [] (v21) at ( \x1, \x2);
		\coordinate [] (v22) at ( \x2, \x2);
		\coordinate [] (v23) at ( \x3, \x2);
		\coordinate [] (v24) at ( \x4, \x2);
		\coordinate [] (v31) at ( \x1, \x3);
		\coordinate [] (v32) at ( \x2, \x3);
		\coordinate [] (v33) at ( \x3, \x3);
		\coordinate [] (v34) at ( \x4, \x3);

		\draw[fill=blue!50!black!20] (v11) -- (v12) -- (v22) -- (v21) -- cycle;
		\draw[fill=red!50!black!20] (v12) -- (v13) -- (v23) -- (v22) -- cycle;
		\draw[fill=red!50!black!20] (v13) -- (v14) -- (v24) -- (v23) -- cycle;
		\draw[fill=blue!50!black!20] (v21) -- (v22) -- (v32) -- (v31) -- cycle;
		\draw[fill=red!50!black!20] (v22) -- (v23) -- (v33) -- (v32) -- cycle;
		\draw[fill=red!50!black!20] (v23) -- (v24) -- (v34) -- (v33) -- cycle;

		\node[label=west:\scriptsize$\dof_1$] at (barycentric cs:v12=0.5,v22=0.5) [] {};
		\node[label=west:\scriptsize$\dof_3$] at (barycentric cs:v22=0.5,v32=0.5) [] {};
		\node[label=south:\scriptsize$\dof_2$] at (barycentric cs:v21=0.5,v22=0.5) [] {};
		\node[label=north:\scriptsize$\dof_4$] at (barycentric cs:v21=0.5,v22=0.5) [] {};
		\node[label=east:\scriptsize$\dof_5$] at (barycentric cs:v12=0.5,v22=0.5) [] {};
		\node[label=east:\scriptsize$\dof_8$] at (barycentric cs:v22=0.5,v32=0.5) [] {};
		\node[label=west:\scriptsize$\dof_6$] at (barycentric cs:v13=0.5,v23=0.5) [] {};
		\node[label=west:\scriptsize$\dof_9$] at (barycentric cs:v23=0.5,v33=0.5) [] {};
		\node[label=south:\scriptsize$\dof_7$] at (barycentric cs:v22=0.5,v23=0.5) [] {};
		\node[label=north:\scriptsize$\dof_{10}$] at (barycentric cs:v22=0.5,v23=0.5) [] {};
		\node[label=east:\scriptsize$\dof_{11}$] at (barycentric cs:v13=0.5,v23=0.5) [] {};
		\node[label=east:\scriptsize$\dof_{13}$] at (barycentric cs:v23=0.5,v33=0.5) [] {};
		\node[label=south:\scriptsize$\dof_{12}$] at (barycentric cs:v23=0.5,v24=0.5) [] {};
		\node[label=north:\scriptsize$\dof_{14}$] at (barycentric cs:v23=0.5,v24=0.5) [] {};

		\end{tikzpicture}

		\caption{Enumeration of dofs ($\dof_i$).}
		\label{fig:relation-table-example-b}
	\end{subfigure}
	\\
	\vspace{5mm}
	\begin{subfigure}[b]{0.45\textwidth}
		\centering
		\begin{tikzpicture}
		[
		cellCentroid/.style={circle,draw=black,text=black,fill=black,thick,minimum size=1mm, inner sep = 0em, outer sep = 0.0em},
		edgeCentroid/.style={rectangle,text=black,thick,minimum size=1mm, inner sep = 0em, outer sep = 0.0em},
		cellLabel/.style={circle,draw=black,text=black,fill=red!50!black!20,thick,minimum size=4mm, inner sep = 0.1em, outer sep = 0.1em},
		edgeLabel/.style={rectangle,draw=black,text=black,fill=blue!50!black!20,thick,minimum size=4mm, inner sep = 0em, outer sep = 0.1em},
		AE1/.style={rectangle,draw=blue!50!black!20,fill=blue!50!black!20,thick,minimum size=4mm, inner sep = 0em, outer sep = 0.0em},
		AE2/.style={rectangle,draw=red!50!black!20,fill=red!50!black!20,thick,minimum size=4mm, inner sep = 0em, outer sep = 0.0em}]
	
		\tikzmath{
			\x1 = 0.0*\unit; \x2 = 0.25*\unit; \x3 = 0.3*\unit; \x4 = 0.55*\unit; \x5 = 0.6*\unit;
			\x6 = 0.85*\unit; \y1 = 0.0*\unit; \y2 = 0.25*\unit; \y3 = 0.5*\unit; \y4 = 0.65*\unit; }

		\coordinate [] (v11) at ( \x1, \y1);
		\coordinate [] (v12) at ( \x2, \y1);
		\coordinate [] (v13) at ( \x3, \y1);
		\coordinate [] (v14) at ( \x4, \y1);
		\coordinate [] (v15) at ( \x5, \y1);
		\coordinate [] (v16) at ( \x6, \y1);
		\coordinate [] (v21) at ( \x1, \y2);
		\coordinate [] (v22) at ( \x2, \y2);
		\coordinate [] (v23) at ( \x3, \y2);
		\coordinate [] (v24) at ( \x4, \y2);
		\coordinate [] (v25) at ( \x5, \y2);
		\coordinate [] (v26) at ( \x6, \y2);
		\coordinate [] (v31) at ( \x1, \y3);
		\coordinate [] (v32) at ( \x2, \y3);
		\coordinate [] (v33) at ( \x3, \y3);
		\coordinate [] (v34) at ( \x4, \y3);
		\coordinate [] (v35) at ( \x5, \y3);
		\coordinate [] (v36) at ( \x6, \y3);
		\coordinate [] (v41) at ( \x1, \y4);
		\coordinate [] (v42) at ( \x2, \y4);
		\coordinate [] (v43) at ( \x3, \y4);
		\coordinate [] (v44) at ( \x4, \y4);
		\coordinate [] (v45) at ( \x5, \y4);
		\coordinate [] (v46) at ( \x6, \y4);

		\draw[fill=blue!50!black!20] (v11) -- (v12) -- (v22) -- (v21) -- cycle;
		\draw[fill=blue!50!black!20] (v21) -- (v22) -- (v32) -- (v31) -- cycle;
		\draw[fill=red!50!black!20] (v13) -- (v14) -- (v24) -- (v23) -- cycle;
		\draw[fill=red!50!black!20] (v23) -- (v24) -- (v34) -- (v33) -- cycle;
		\draw[fill=red!50!black!20] (v15) -- (v16) -- (v26) -- (v25) -- cycle;
		\draw[fill=red!50!black!20] (v25) -- (v26) -- (v36) -- (v35) -- cycle;

		\draw[draw=red!50!black,very thick,dashed] (barycentric cs:v12=0.5,v13=0.5) -- (barycentric cs:v42=0.5,v43=0.5);
		\draw[draw=red!50!black,very thick,dashed] (barycentric cs:v14=0.5,v15=0.5) -- (barycentric cs:v44=0.5,v45=0.5);

		\node[] at (barycentric cs:v11=0.5,v22=0.5) [] {$\elem_1$};
		\node[] at (barycentric cs:v21=0.5,v32=0.5) [] {$\elem_2$};
		\node[] at (barycentric cs:v13=0.5,v24=0.5) [] {$\elem_3$};
		\node[] at (barycentric cs:v23=0.5,v34=0.5) [] {$\elem_4$};
		\node[] at (barycentric cs:v15=0.5,v26=0.5) [] {$\elem_5$};
		\node[] at (barycentric cs:v25=0.5,v36=0.5) [] {$\elem_6$};
		
		\node[] at (barycentric cs:v31=0.5,v42=0.5) [] {\color{red!50!black}$C_1$};
		\node[] at (barycentric cs:v33=0.5,v44=0.5) [] {\color{red!50!black}$C_2$};
		\node[] at (barycentric cs:v35=0.5,v46=0.5) [] {\color{red!50!black}$C_3$};

		\end{tikzpicture}

		\caption{Initial parallel distribution of elements in 3 cores.}
		\label{fig:relation-table-example-c}
	\end{subfigure}
	\hfill
	\begin{subfigure}[b]{0.45\textwidth}
		\centering
		\begin{tikzpicture}
		[
		cellCentroid/.style={circle,draw=black,text=black,fill=black,thick,minimum size=1mm, inner sep = 0em, outer sep = 0.0em},
		edgeCentroid/.style={rectangle,text=black,thick,minimum size=1mm, inner sep = 0em, outer sep = 0.0em},
		cellLabel/.style={circle,draw=black,text=black,fill=red!50!black!20,thick,minimum size=4mm, inner sep = 0.1em, outer sep = 0.1em},
		edgeLabel/.style={rectangle,draw=black,text=black,fill=blue!50!black!20,thick,minimum size=4mm, inner sep = 0em, outer sep = 0.1em},
		AE1/.style={rectangle,draw=blue!50!black!20,fill=blue!50!black!20,thick,minimum size=4mm, inner sep = 0em, outer sep = 0.0em},
		AE2/.style={rectangle,draw=red!50!black!20,fill=red!50!black!20,thick,minimum size=4mm, inner sep = 0em, outer sep = 0.0em}]
		
		\tikzmath{
			\x1 = 0.0*\unit; \x2 = 0.25*\unit; \x3 = 0.3*\unit; \x4 = 0.55*\unit; \x5 = 0.8*\unit;
			\x6 = 0.85*\unit; \x7 = 0.95*\unit; \y1 = 0.0*\unit; \y2 = 0.25*\unit; \y3 = 0.5*\unit; \y4 = 0.65*\unit; }

		\coordinate [] (v11) at ( \x1, \y1);
		\coordinate [] (v12) at ( \x2, \y1);
		\coordinate [] (v13) at ( \x3, \y1);
		\coordinate [] (v14) at ( \x4, \y1);
		\coordinate [] (v15) at ( \x5, \y1);
		\coordinate [] (v16) at ( \x6, \y1);
		\coordinate [] (v21) at ( \x1, \y2);
		\coordinate [] (v22) at ( \x2, \y2);
		\coordinate [] (v23) at ( \x3, \y2);
		\coordinate [] (v24) at ( \x4, \y2);
		\coordinate [] (v25) at ( \x5, \y2);
		\coordinate [] (v26) at ( \x6, \y2);
		\coordinate [] (v31) at ( \x1, \y3);
		\coordinate [] (v32) at ( \x2, \y3);
		\coordinate [] (v33) at ( \x3, \y3);
		\coordinate [] (v34) at ( \x4, \y3);
		\coordinate [] (v35) at ( \x5, \y3);
		\coordinate [] (v36) at ( \x6, \y3);
		\coordinate [] (v41) at ( \x1, \y4);
		\coordinate [] (v42) at ( \x2, \y4);
		\coordinate [] (v43) at ( \x3, \y4);
		\coordinate [] (v44) at ( \x4, \y4);
		\coordinate [] (v45) at ( \x5, \y4);
		\coordinate [] (v46) at ( \x6, \y4);
		\coordinate [] (v37) at ( \x7, \y3);
		\coordinate [] (v47) at ( \x7, \y4);

		\draw[fill=blue!50!black!20] (v11) -- (v12) -- (v22) -- (v21) -- cycle;
		\draw[fill=blue!50!black!20] (v21) -- (v22) -- (v32) -- (v31) -- cycle;
		\draw[fill=red!50!black!20] (v13) -- (v14) -- (v24) -- (v23) -- cycle;
		\draw[fill=red!50!black!20] (v23) -- (v24) -- (v34) -- (v33) -- cycle;
		\draw[fill=red!50!black!20] (v14) -- (v15) -- (v25) -- (v24) -- cycle;
		\draw[fill=red!50!black!20] (v24) -- (v25) -- (v35) -- (v34) -- cycle;

		\draw[draw=red!50!black,very thick,dashed] (barycentric cs:v12=0.5,v13=0.5) -- (barycentric cs:v42=0.5,v43=0.5);
		\draw[draw=red!50!black,very thick,dashed] (barycentric cs:v15=0.5,v16=0.5) -- (barycentric cs:v45=0.5,v46=0.5);

		\node[] at (barycentric cs:v11=0.5,v22=0.5) [] {$\elem_1$};
		\node[] at (barycentric cs:v21=0.5,v32=0.5) [] {$\elem_2$};
		\node[] at (barycentric cs:v13=0.5,v24=0.5) [] {$\elem_3$};
		\node[] at (barycentric cs:v23=0.5,v34=0.5) [] {$\elem_4$};
		\node[] at (barycentric cs:v14=0.5,v25=0.5) [] {$\elem_5$};
		\node[] at (barycentric cs:v24=0.5,v35=0.5) [] {$\elem_6$};
		
		\node[] at (barycentric cs:v31=0.5,v42=0.5) [] {\color{red!50!black}$C_1$};
		\node[] at (barycentric cs:v33=0.5,v45=0.5) [] {\color{red!50!black}$C_2$};
		\node[] at (barycentric cs:v37=0.5,v47=0.5) [] {\color{red!50!black}$C_3$};

		\end{tikzpicture}

		\caption{Parallel distribution of elements after redistribution.}
		\label{fig:relation-table-example-d}
	\end{subfigure}
	\\
	\vspace{5mm}
	\begin{subfigure}[b]{0.48\textwidth}
		\centering
		\begin{tikzpicture}
		\matrix[
			matrix of nodes,
			text height=1.6ex,
			text depth=0.0ex,
			text width=1.6ex,
			align=center,
			nodes={draw=black!5}, 
			nodes in empty cells,
			left delimiter={[},
			right delimiter={]}
		] at (0,0) (A){
			&&&&&&&&&&&&&\\
			&&&&&&&&&&&&&\\
			&&&&&&&&&&&&&\\
			&&&&&&&&&&&&&\\
			&&&&&&&&&&&&&\\
			&&&&&&&&&&&&&\\
		};
		\node[] at (A-1-1){1};
		\node[] at (A-1-2){1};
		\node[] at (A-2-3){1};
		\node[] at (A-2-4){1};
		\node[] at (A-3-5){1};
		\node[] at (A-3-6){1};
		\node[] at (A-3-7){1};
		\node[] at (A-4-8){1};
		\node[] at (A-4-9){1};
		\node[] at (A-4-10){1};
		\node[] at (A-5-11){1};
		\node[] at (A-5-12){1};
		\node[] at (A-6-13){1};
		\node[] at (A-6-14){1};
		\foreach \i in {1,...,14}
			{
			\node[anchor=south] at (A-1-\i.north){$\dof_{\i}$};
			}
		\foreach \i in {1,...,6}
		{
			\node[anchor=east] at (A-\i-1.west){$\elem_\i{\color{white}\hspace{3mm}}$};
		}    
		\end{tikzpicture}
		\caption{``element\_dof" relation.}
		\label{fig:element-dof}
	\end{subfigure}
	\hfill
	\begin{subfigure}[b]{0.48\textwidth}
		\centering
		\begin{tikzpicture}
		\matrix[
			matrix of nodes,
			text height=1.6ex,
			text depth=0.0ex,
			text width=1.6ex,
			align=center,
			nodes={draw=black!5}, 
			nodes in empty cells,
			left delimiter={[},
			right delimiter={]}
		] at (0,0) (A){
			&&&&&&&&&&&&&\\
			&&&&&&&&&&&&&\\
			&&&&&&&&&&&&&\\
			&&&&&&&&&&&&&\\
			&&&&&&&&&&&&&\\
			&&&&&&&&&&&&&\\
			&&&&&&&&&&&&&\\
		};
		\node[] at (A-1-1){1};
		\node[] at (A-3-2){1};
		\node[] at (A-2-3){1};
		\node[] at (A-3-4){1};
		\node[] at (A-1-5){1};
		\node[] at (A-4-6){1};
		\node[] at (A-6-7){1};
		\node[] at (A-2-8){1};
		\node[] at (A-5-9){1};
		\node[] at (A-6-10){1};
		\node[] at (A-4-11){1};
		\node[] at (A-7-12){1};
		\node[] at (A-5-13){1};
		\node[] at (A-7-14){1};
		\foreach \i in {1,...,14}
			{
			\node[anchor=south] at (A-1-\i.north){$\dof_{\i}$};
			}
		\foreach \i in {1,...,7}
		{
			\node[anchor=east] at (A-\i-1.west){$\truedof_\i{\color{white}\hspace{3mm}}$};
		}    
		\end{tikzpicture}
		\caption{(``dof\_truedof"$)^T$ relation.}
		\label{fig:truedof-dof}
	\end{subfigure}
	\\
	\vspace{5mm}
	\begin{subfigure}[b]{0.45\textwidth}
		\centering
		\begin{tikzpicture}
		\matrix[
			matrix of nodes,
			text height=2ex,
			text depth=0.0ex,
			text width=2ex,
			align=center,
			nodes={draw=black!5}, 
			nodes in empty cells,
			left delimiter={[},
			right delimiter={]}
		] at (0,0) (A){
			&&&&&\\
			&&&&&\\
			&&&&&\\
		};
		\node[] at (A-1-1){1};
		\node[] at (A-1-2){1};
		\node[] at (A-2-3){1};
		\node[] at (A-2-4){1};
		\node[] at (A-3-5){1};
		\node[] at (A-3-6){1};

		\foreach \i in {1,...,6}
			{
			\node[anchor=south] at (A-1-\i.north){$\elem_\i$};
			}
		\foreach \i in {1,...,3}
		{
			\node[anchor=east] at (A-\i-1.west){$C_\i{\color{white}\hspace{3mm}}$};
		}
		\end{tikzpicture}
		\caption{``core\_element" relation in Algorithm~\ref{algorithm: parallel agglomeration}.}
		\label{fig:core-element}
	\end{subfigure}
	\hfill
	\begin{subfigure}[b]{0.45\textwidth}
		\centering
		\begin{tikzpicture}
		\matrix[
			matrix of nodes,
			text height=2ex,
			text depth=0.0ex,
			text width=2ex,
			align=center,
			nodes={draw=black!5}, 
			nodes in empty cells,
			left delimiter={[},
			right delimiter={]}
		] at (0,0) (A){
			&&\\
			&&\\
			&&\\
		};
		\node[] at (A-1-1){1};
		\node[] at (A-2-2){1};
		\node[] at (A-2-3){1};

		\foreach \i in {1,...,3}
			{
			\node[anchor=south] at (A-1-\i.north){$C_\i$};
			}
		\foreach \i in {1,...,3}
		{
			\node[anchor=east] at (A-\i-1.west){$C_\i{\color{white}\hspace{3mm}}$};
		} 
		\end{tikzpicture}
		\caption{``Core\_core" relation in Algorithm~\ref{algorithm: parallel agglomeration}.}
		\label{fig:Core-core}
	\end{subfigure}
	\\
	\vspace{5mm}
	\begin{subfigure}[b]{0.45\textwidth}
		\centering
		\begin{tikzpicture}
		\matrix[
			matrix of nodes,
			text height=2ex,
			text depth=0.0ex,
			text width=2ex,
			align=center,
			nodes={draw=black!5}, 
			nodes in empty cells,
			left delimiter={[},
			right delimiter={]}
		] at (0,0) (A){
			&&&&&\\
			&&&&&\\
			&&&&&\\
		};
		\node[] at (A-1-1){1};
		\node[] at (A-1-2){1};
		\node[] at (A-2-3){1};
		\node[] at (A-2-4){1};
		\node[] at (A-2-5){1};
		\node[] at (A-2-6){1};

		\foreach \i in {1,...,6}
			{
			\node[anchor=south] at (A-1-\i.north){$\elem_\i$};
			}
		\foreach \i in {1,...,3}
		{
			\node[anchor=east] at (A-\i-1.west){$C_\i{\color{white}\hspace{3mm}}$};
		}
		\end{tikzpicture}
		\caption{``Core\_element" = ``Core\_core" $\times$ ``core\_element".}
		\label{fig:Core-element}
	\end{subfigure}
	\hfill
	\begin{subfigure}[b]{0.45\textwidth}
		\centering
		\begin{tikzpicture}
		\matrix[
			matrix of nodes,
			text height=2ex,
			text depth=0.0ex,
			text width=2ex,
			align=center,
			nodes={draw=black!5}, 
			nodes in empty cells,
			left delimiter={[},
			right delimiter={]}
		] at (0,0) (A){
			&&&&&\\
			&&&&&\\
		};
		\node[] at (A-1-1){1};
		\node[] at (A-1-2){1};
		\node[] at (A-2-3){1};
		\node[] at (A-2-4){1};
		\node[] at (A-2-5){1};
		\node[] at (A-2-6){1};

		\foreach \i in {1,...,6}
			{
			\node[anchor=south] at (A-1-\i.north){$\elem_\i$};
			}
		\foreach \i in {1,...,2}
		{
			\node[anchor=east] at (A-\i-1.west){$AE_\i{\color{white}\hspace{3mm}}$};
		}
		\end{tikzpicture}
		\caption{``AE\_element" relation.}
		\label{fig:AE-element}
	\end{subfigure}

	\caption{An example of mesh elements, the associated Raviart-Thomas dofs, true dofs, AEs, and relevant relation tables.}
	\label{fig:relation-table-example}
\end{figure}

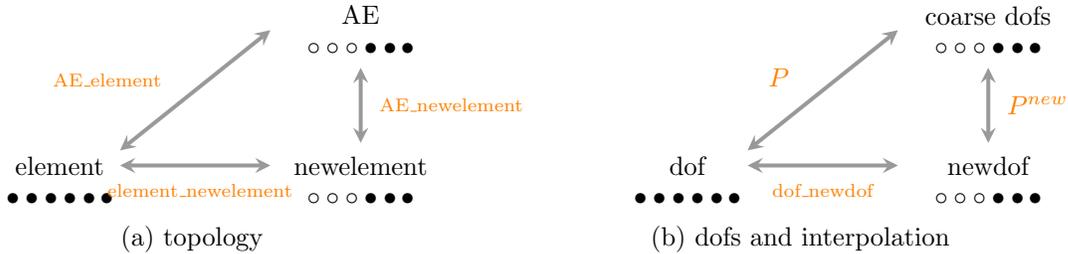
\begin{figure}[ht]
\centering
\begin{tikzpicture}

\coordinate [label=] (fine) at (0, 0);
\coordinate [label=] (newfine) at (4, 0);
\coordinate [label=] (coarse) at (4, 2);

\node[label={south:$\bullet\bullet\bullet\bullet\bullet\,\bullet$}] at (fine) [] {element};
\node[label={south:$\circ\circ\circ\bullet\bullet\,\bullet$}] at (newfine) [] {newelement};
\node[label={south:$\circ\circ\circ\bullet\bullet\,\bullet$}] at (coarse) [] {AE};

\node[label={east:{\color{orange}\scriptsize AE\_newelement}}] at (barycentric cs:coarse=0.4,newfine=0.6) [] {};
\node[label={north west:{\color{orange}\scriptsize AE\_element}}] at (barycentric cs:coarse=0.4,fine=0.6) [] {};
\node[label={south:{\color{orange}\scriptsize element\_newelement}}] at (barycentric cs:newfine=0.48,fine=0.55) [] {};

\draw[black!40!, line width= 1.5pt, stealth-stealth] (barycentric cs:fine=0.8,newfine=0.2) -- (barycentric cs:fine=0.3,newfine=0.7);
\draw[black!40!, line width= 1.5pt, stealth-stealth] (barycentric cs:coarse=0.65,newfine=0.35) -- (barycentric cs:coarse=0.15,newfine=0.85);
\draw[black!40!, line width= 1.5pt, stealth-stealth] (2.8, 1.8) -- (0.8, 0.2);

\end{tikzpicture}
\hspace{1cm}
\begin{tikzpicture}

\coordinate [label=] (fine) at (0, 0);
\coordinate [label=] (newfine) at (4, 0);
\coordinate [label=] (coarse) at (4, 2);

\node[label={south:$\bullet\bullet\bullet\bullet\bullet\,\bullet$}] at (fine) [] {dof};
\node[label={south:$\circ\circ\circ\bullet\bullet\,\bullet$}] at (newfine) [] {newdof};
\node[label={south:$\circ\circ\circ\bullet\bullet\,\bullet$}] at (coarse) [] {coarse dofs};

\node[label={east:{\color{orange}$P^{new}$}}] at (barycentric cs:coarse=0.4,newfine=0.6) [] {};
\node[label={north west:{\color{orange}$P$}}] at (barycentric cs:coarse=0.4,fine=0.6) [] {};
\node[label={south:{\color{orange}\scriptsize dof\_newdof}}] at (barycentric cs:newfine=0.45,fine=0.55) [] {};

\draw[black!40!, line width= 1.5pt, stealth-stealth] (barycentric cs:fine=0.8,newfine=0.2) -- (barycentric cs:fine=0.3,newfine=0.7);
\draw[black!40!, line width= 1.5pt, stealth-stealth] (barycentric cs:coarse=0.65,newfine=0.35) -- (barycentric cs:coarse=0.15,newfine=0.85);
\draw[black!40!, line width= 1.5pt, stealth-stealth] (2.8, 1.8) -- (0.8, 0.2);

\end{tikzpicture}\\
(a) topology \hspace*{140pt} (b) dofs \updates{and interpolation} \\
\caption{An Illustration of (a) how ``element\_newelement", ``AE\_newelement", and ``AE\_element" connect element, ``newelement", and AE, as in Algorithm~\ref{algorithm: parallel agglomeration}; and (b) how ``dof\_newdof"$ = ($``newdof\_dof"$)^T$, $P^{new}$, and $P$ connect dof, ``newdof", and coarse dofs, as in Algorithm~\ref{algorithm: data redistribution}. Here, $\bullet$ and $\circ$ represent active and inactive cores respectively.}
\label{fig:newP}
\end{figure}

With data redistribution for coarse levels, the agglomerated elements -- from the relation ``AE\_element" -- may cover elements from different cores; that is, elements may be redistributed to a different (in our case, smaller) set of cores. Since the coarsening algorithm in \cite{deRhamAMGe} requires each AE to be path-connected in the physical domain, we need to make sure the elements in each core after redistribution are physically connected. To this end, we make use of the relation ``core\_element", which lists the active cores and the elements associated with them \updates{(see Figures~\ref{fig:relation-table-example-c} and \ref{fig:core-element} for an example)}, to guide the redistribution. In what follows, we describe how the new relations, which enable redistribution, are formed. This is accomplished with regards to the topology in Algorithm~\ref{algorithm: parallel agglomeration} and with regards to the dofs in Algorithm~\ref{algorithm: data redistribution}. We refer the reader to Figure~\ref{fig:newP} for an illustration of these algorithms.

\begin{algorithm}[Agglomeration in parallel with redistribution]\label{algorithm: parallel agglomeration}\hfill

\begin{itemize}
\item We form the new relation:\\
``core\_core" $=$ ``core\_element" $\times $ ``element\_element" $\times ($``core\_element"$)^T$.
\item We use the relation (graph) ``core\_core" as an input for a graph partitioner and create a new relation ``Core\_core" \updates{(see Figure~\ref{fig:Core-core})}, which aggregates the currently active  cores to a smaller set of ``Cores" (that can be a subset of the current ones). 
\item Then, we form ``Core\_element" $=$ ``Core\_core" $\times$ ``core\_element" \updates{(see Figure~\ref{fig:Core-element})}.
\item Based on ``Core\_element" we make a copy of each element related to a ``Core" denoted by ``newelement". That is, we create the relation 
``newelement\_element", which effectively redistributes the set of elements from the current configuration of cores to the new configuration of ``Cores". 
\item Finally, on each ``Core", we use a local graph partitioner to create ``AE\_newelement" using as input the relation 
``newelement\_newelement" $=$ ``newelement\_element" $\times$ ``element\_element"$\times ($``newelement\_element"$)^T$.
\item The desired output is ``AE\_element" $=$``AE\_newelement"$\times$``newelement\_element" in addition to ``newelement\_element" \updates{(see Figure~\ref{fig:relation-table-example-d},\subref{fig:AE-element})}. 
\end{itemize}
\end{algorithm}

One user-specified parameter for Algorithm~\ref{algorithm: parallel agglomeration} is the core coarsening factor $\beta_c$, which is the ratio of number of cores to number of ``Cores". Given $\beta_c$, number of ``Cores" is calculated as
\begin{equation}
\text{number of ``Cores"} = \lceil \;  \text{number of cores} \,/\, \beta_c \; \rceil,
\label{eq:core-coarsening-factor}
\end{equation}
where $\lceil x\rceil$ is the smallest integer that is greater than or equal to $x$. \chak{A simple example showing the distributions of elements before and after redistribution, and the corresponding relations in Algorithm~\ref{algorithm: parallel agglomeration} can be found in Figure~\ref{fig:relation-table-example}.}

Next, we want to build coarse spaces associated with the new AEs\updates{, addressing the fact that the AEs can cover elements from multiple cores}. Since the  construction of $\{P_T\}$  for each agglomerate $T$ are local, the option we have chosen is to redistribute \updates{(i.e., relocate)} the \updates{fine-level} elements that form $T$  to the core that $T$ is assigned to. 

\begin{algorithm}[Redistribution algorithm]\label{algorithm: data redistribution} \hfill

We assume that we have $A_{\diag}$ in terms of $A_\tau$, all in terms of (decoupled) \updates{fine-level} dofs. We also have the ``dof\_truedof" and  ``element\_dof" relations  \updates{(see Figure~\ref{fig:element-dof}-\subref{fig:truedof-dof})}. 

The  data redistribution in this algorithm is guided by the two parallel relations ``AE\_element" and ``newelement\_element" constructed, for example, 
using Algorithm~\ref{algorithm: parallel agglomeration}.

We form the following relations using  parallel matrix computations:
\begin{itemize}
\item Form the composite relation ``AE\_dof" $=$``AE\_element" $\times$ ``element\_dof".
\item From the relation ``AE\_dof", we form a new relation: ``newdof\_dof", where each ``newdof" is a copy of the dof that belongs to its \updates{respective} unique\footnote{\updates{Recall that each (decoupled) fine-level dof is uniquely assigned to a fine-level element and consequently it is uniquelly assigned to the AE that contains its element in accordance with ``AE\_element". Generally, ``AE\_dof" represents a relation that can cross core boundaries, while the copies in ``newdof" address that so that ``AE\_newdof" relates only AEs and ``newdofs" located on the same cores.}} AE but assigned to the core that owns that AE. This effectively redistributes the dofs from one configuration of cores to the new configuration of cores \updates{induced by} the distribution of the AEs.

\item We form the relation\footnote{\updates{Clearly, ``newelement\_newdof" relates only ``newelements" and ``newdofs" located on the same cores.}} ``newelement\_newdof" = ``newelement\_element" $\times$ ``element\_dof" $\times$ $(``newdof\_dof")^T$. 
\item \updates{Element matrices are redistributed by computing} numerically the product ``newdof\_dof" $\times A_{\diag} \times ($``newdof\_dof"$)^T$ which is the redistribution of $A_{\diag}$, denoted by $A^{new}_{\diag}$.
\item We need also to connect the coarse level at the new core configuration to the fine level at the previous core configuration. This is achieved by the product
$P$ = $(``newdof\_dof")^T \times P^{new}$, where $P^{new}$ is an interpolation matrix constructed at the new core configuration \updates{in terms of ``newdof", e.g., using the procedure discussed in Section~\ref{section: parallel finite elements}}. That is, $P$ relates the coarse dofs from the new core configuration with the (fine) dofs from the previous core configuration. See Figure~\ref{fig:newP} for an illustration.

\item Finally, we can create the ``newtruedof" from the products ``newdof\_truedof" $=$ ``newdof\_dof"$\times$``dof\_truedof"
 and ``newdof\_newdof" $=$ ``newdof\_truedof" $\times($``newdof\_truedof"$)^T$. From the latter one\footnote{\updates{In ``newdof\_newdof", the ``newdofs" are related through the true dofs, in the sense that two ``newdofs" are related if and only if they are decoupled versions of the same true dof. Thus, the ``newdofs" form equivalence classes, where each class represents a true dof. In practice, we select one of the mutually-related ``newdofs" to be a representative of the class and hence of the respective true dof. This effectively produces the ``newtruedof" redistribution of the original ``truedof", where both conceptually represent the same global true dofs but in different configurations.}}, we select one of the ``newdofs" to be an actual ``newtruedof" and connect it  with the remaining ``newdofs", i.e., we create  the relation ``newtruedof\_newdof".  
The relation ``newtruedof \_truedof" $=$ ``newtruedof\_newdof" $\times$ ``newdof\_truedof" effectively redistributes the true dofs. 
\end{itemize} 
\end{algorithm}

\begin{figure}
	\centering
	\includegraphics[width=.3\textwidth]{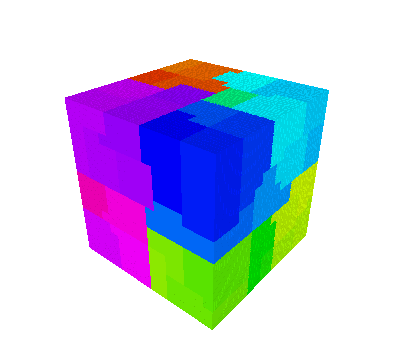}
	\includegraphics[width=.3\textwidth]{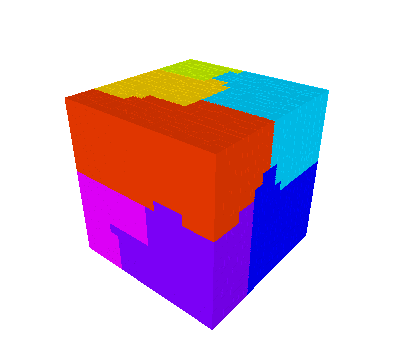}
	\includegraphics[width=.3\textwidth]{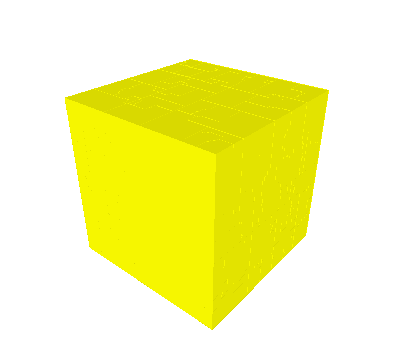}\\
	(a) $\text{nc}=64$ \hspace*{85pt}(b) $\text{nc}=8$\hspace*{85pt}(c) $\text{nc}=1$
	\caption{Schematic of data redistribution on three levels for a unit cube spatial domain. Each color denotes an independent computational process appointed to a respective core, with (a) $\text{nc}=64$, (b) $\text{nc}=8$, and (c) $\text{nc}=1$ as the number of cores to which the mesh is redistributed for that level.
	 }\label{fig:cubes}
\end{figure}

Once we have all element and global matrix data and relations redistributed with respect to the new core configuration, we apply the AMGe algorithm as implemented in the case of \updates{a steady} configuration of cores  since now the computations can be performed locally with respect to the cores from the new configuration\updates{, and the AMGe algorithm is agnostic about the existence of other core configurations}. Note that some cores can become idle at this point. 
An important point to note is that  after the redistribution,  the data from the most recent \updates{level} of AMGe on \updates{the} previous core configuration is no longer needed\updates{, since it is only used temporarily for the construction procedure}. A visual of the redistribution applied to multiple levels of a hierarchy is provided in Figure~\ref{fig:cubes}.

\begin{remark}[Comments on computational time]\label{sec:offlinecost}
As this redistribution functionality is integrated into the AMGe hierarchy build stage (for both topology and finite element spaces) through the formation of relations in Algorithm~\ref{algorithm: parallel agglomeration} and Algorithm \ref{algorithm: data redistribution}, the total walltime to perform redistribution is treated as a one-time, offline cost. We note, there are some \updates{small} increases in online computational time (\updates{during} simulation), such as \updates{increases} associated with composite calculations to evaluate the bilinear form, as well as those which are consequences of redistribution, i.e., fewer cores means an increase in local problem size which increases the wall time. These online costs, which are dependent on the particular redistribution configuration, are included in the walltime per simulation when evaluating computational performance.
\end{remark}

\section{MC and MLMC methods for estimating statistical moments of QoI involving solutions of PDEs with random coefficients}\label{section: MLMC application}
We consider the following model problem, where for a given  permeability field $k$ defined on a domain $D$, the pressure $p\in L^2(D)$ and velocity field $\bm q\in H(\div , D)$ are solved via the system of PDEs:
\begin{equation}\label{eq:darcy}
 			\left\{
  				\begin{array}{lcl}
    				k^{-1}\bm q + \nabla p = \bm f, &\ \ & \text{in }D, \vspace{1mm}\\
    				\nabla \cdot  \bm q = 0, &\ \ &\text{in }D , %
			  \end{array}
			\right.
\end{equation}	 
with $p=p_D$ on $\Gamma_D$, enforced by $\bm f$, and $\bm q \cdot \bm n = 0$ on $\Gamma_N$, where $\Gamma_D$ and $\Gamma_N$ are non-overlapping partitions of $\partial D$. Due to a lack of information regarding the permeability field, $k$ is treated as a source of uncertainty. Given the probability space $(\Omega, \mathcal{F}, P)$, with sample space $\Omega$, the permeability field is defined as $k:=\exp(u(\bm x, \omega))$ with random realizations $u$ from $\{u(\bm x, \omega) \in L^2(D):  \bm x\in D, \omega \in \Omega\}$ distributed according to a zero-mean Gaussian distribution. %

Consequently, solutions of velocity, defined $\bm q:=\bm q(\bm x, \omega)$ are random, and the task is to estimate the effect of the random coefficient on solutions of interest. In particular, we seek to estimate statistical moments of the scalar QoI (quantity of interest) $Q(\omega)$ representing the flux across the outflow boundary $\Gamma_\text{out} \subset \Gamma_D$, and defined as
\begin{equation}\label{eq:Q}
Q(\omega) = \frac{1}{|\Gamma_\text{out}|}\int_{\Gamma_\text{out}} \bm q (\cdot, \omega)\cdot \bm n \ dS,
\end{equation}
using sampling-based uncertainty quantification (UQ) strategies. This is accomplished by performing repeated PDE simulations with various realizations of the permeability field. 
In this work, we utilize numerical approximations of $\bm q$ (and hence $Q$) generated by solving the mixed PDE in (\ref{eq:darcy}). Given a mesh \hillary{$\mathcal{T}_M$ with $M$ elements (of mesh size $h$)}, the finite dimensional $L^2$- and $H(\div)$-conforming spaces are defined as $\Theta_M$ and $\R_M$; in this work we consider $\Theta _M$ to be the space of piecewise constant elements, and $\R_M$ to be the lowest order Raviart-Thomas elements, with resulting discrete approximations $\bm q_M$ and $Q_M$.

A straight-forward approach to estimate statistical moments of $Q$, and in particular $\mathbb{E}[Q]$, is via a Monte Carlo sampling approach, whereby independently sampled discrete approximations $Q_M^{(i)}$ are averaged as in:
\begin{equation}\label{eq:slmc}
\hat{Q}_M^{MC}= \frac{1}{N}\sum\limits_{i=1}^N Q_M^{(i)}.
\end{equation}
\hillary{The resulting mean square error (MSE) may be decomposed into the statistical error and the square of the discretization error (bias):}
\[
\mathbb{E}[(\hat{Q}_M^{MC}-Q)^2]= N^{-1}\Var(Q_M)+ (\mathbb{E}[{Q}_M-Q])^2.
\]
\hillary{Assuming the discretization error converges as $M^{-\alpha}$, where $\alpha$ denotes the order of convergence, then an MSE bound $\varepsilon^2$ is assumed to be met if $N\geq  2\varepsilon^{-2}\Var(Q_M)$ and $ M\geq (\varepsilon/\sqrt{2})^{-1/\alpha}$, where the error is balanced between the two terms with.}

As the number of fine simulations depends on the variance $\Var(Q_M)$, standard Monte Carlo is not suited for large-scale applications where generating a large number of simulations is infeasible. To alleviate this computational burden, we consider Multilevel Monte Carlo.

\subsection{Multilevel Monte Carlo}
\label{ssec:MLMC}
As a variance reduction approach amenable to large-scale UQ, multilevel Monte Carlo~\cite{Heinrich01,Giles08, Giles13, barth2011multi, Cliffe11} is applied to estimate statistical moments of $Q$. To do this, the fine grid mesh $\T_0:=\T_M$ is recursively agglomerated to $L$ coarser levels $\{\T_\ell\}_{\ell=1}^L$, \hillary{with decreasing number of elements $M:=M_0 > M_1>\hdots > M_L$ (equivalently, increasing mesh size),} such that $\T_L$ is the coarsest. For a given $\T_\ell$, $\Theta_\ell$ and $\R_\ell$ are the associated pair of lowest-order Raviart-Thomas finite element spaces.

On each level, discrete solutions of the
Darcy velocity and the QoI for level $\ell=0,1,\hdots L$ are denoted $\bm q_\ell$ and $Q_\ell$, respectively. Using a telescoping sum, the expected value of $Q_0$ may be rewritten as
\begin{equation}\label{eq:QML-unused}
\mathbb{E} [Q_0] = \mathbb{E}[Q_L]  + \sum\limits_{\ell=0}^{L-1} \mathbb{E}[Q_\ell  -Q_{\ell+1}].
\end{equation}
In what follows, we define $Y_L=Q_L$ and $Y_\ell = Q_\ell - Q_{\ell+1}$ for $\ell = 0,1,\hdots ,L-1$. To estimate the terms of the telescoping sum, a Monte Carlo (MC) estimate is applied to each expectation. In other words, 
\begin{equation}\label{eq:mlmc}
\hat{Q}_0^{ML}:= \sum\limits_{\ell=0}^{L} \hat{Y}_\ell^{MC},
\end{equation}
such that coarsest term $\hat{Y}_L^{MC}$ is estimated as in (\ref{eq:slmc}), and for $\ell=0,
\hdots, L-1$,
\begin{equation}\label{eq:QML}
 \hat{Y}_\ell^{MC}:= \frac{1}{N_\ell}\sum\limits_{i=1}^{N_\ell} \left(Q_\ell^{(i)}-Q_{\ell+1}^{(i)}\right)
\end{equation}
is a Monte Carlo estimate utilizing the error between the fine and coarse samples $Q_\ell^{(i)}$ and $Q_{\ell+1}^{(i)}$ estimated from matching fine and coarse permeability fields. %

The efficiency gain of MLMC (in comparison to MC), comes from the ability to distribute the computational work amongst the various MC estimates in (\ref{eq:mlmc}). Define $\C _\ell$ as the computational time to sample $Q_L^{(i)}$ when $\ell=L$ and the combined computational time to sample $Q_\ell^{(i)}$ and $Q_{\ell+1}^{(i)}$ for $\ell \in \{0,1,\hdots , L-1 \}$. Then it can be shown that
the optimal number of samples on each level is given by the formula
\begin{equation}\label{eq:nl}
    N_\ell = 2\varepsilon^{-2} \sqrt{\frac{\Var(Y_\ell)}{\C_\ell}}\sum\limits^{L}_{k=0} \sqrt{\C_k \Var(Y_k)}
\end{equation}
{where $\varepsilon^2$ is the desired MSE of the estimator.}
The above choice of number of samples balances the computational work across all levels in the hierarchy in order to minimize the total computational work~\cite{Giles08}.


\subsection{Cost improvement with Coarsening}\label{sec:mlmctheory}

\hillary{The MLMC theory of~\cite{Giles08,Cliffe11} seeks to identify the number of level refinements needed to reach a desired MSE tolerance. 
In what follows, we denote $nc_\ell$ as the number of cores utilized on level $\ell$.
Suppose there are positive constants $\alpha,\beta,\gamma, c_1,c_2,c_3>0$ such that $\alpha \geq \frac{1}{2} \min (\beta,\gamma)$ and
\begin{enumerate}
\item $\vert E[Q_{\ell}-Q]\vert \leq c_1 M_\ell^{-\alpha}$
\item $V(Y_\ell)  \leq c_2  M^{-\beta}_\ell$
\item $nc_\ell C_\ell \leq c_3  M_\ell^{\gamma}$
\end{enumerate}
Then for any sufficiently small $\varepsilon >0$ there exists a fine level with corresponding $M=M_0$ such that
$$ E[(\hat{Q}^{ML}-Q)^2] < \varepsilon^2$$
and the cost is estimated 
\[
C(\hat{Q}^{ML}) \lesssim
 \left\{\begin{array}{lr}
\varepsilon^{-2}, & \text{ if } \beta > \gamma \\
\varepsilon^{-2}(\log \varepsilon)^2, & \text{ if } \ \beta =\gamma \\
\varepsilon^{-2-(\gamma-\beta)/\alpha}, & \text{ if } \ \beta < \gamma \\
\end{array}
\right.
\]
with ``$\lesssim$'' indicating that the total cost will scale proportionally with changes in $\varepsilon$.}

\hillary{While this is important to verify for the MLMC results, in this work, we consider the opposite problem: for a fixed fine level (that meets the desired discretization error), when should a coarser level be incorporated into the MLMC hierarchy to minimize the cost? As discussed previously, the extent to which the hierarchy may be coarsened is constrained by the number of processors utilized. Given that we no longer have that restriction -- with the introduction of level-dependent data distribution -- we seek to identify the configuration that theoretically improves the cost. To that end, we build our cost estimates from the traditional MLMC theory assumptions (see, e.g., ~\cite{Cliffe11}).}

\hillary{The total cost of MLMC with $L+1$ levels can be calculated by
\begin{equation}\label{eq:MLMCcost}
C^\text{total}_{L} = \sum\limits_{\ell=0}^L N_\ell \mathcal{C}_\ell.
\end{equation}
Using the definition of $N_\ell$ (from \eqref{eq:nl}) and applying the relationship $V(Y_\ell)  = c_2  M^{-\beta}_\ell$ it can be shown that the cost can be approximated as
\begin{equation}\nonumber
\hat{C}^\text{total}_{L} =  2\varepsilon^{-2}\left( \sqrt{V(Q_L)C_L}+ \sum\limits_{\ell=0}^{L-1} \sqrt{c_2  M^{-\beta}_\ell\mathcal{C}_\ell} \right)^2.
\end{equation}}

\hillary{In practice, given a partial hierarchy of levels $\ell=0,1,\hdots, K$, with $K\geq 2$, we can estimate if there will be cost improvement by increasing $K$. This requires estimates of $V(Q_\ell)$ and $V(Y_\ell)$ from the partial hierarchy, as well as predictions for the cost on each level in proposed full hierarchy with levels $\ell=0,1,\hdots, L$. Here we define the {\it Scaled Root Cost for $L+1$ levels as as}
\begin{equation}\label{eq:scaledrootcost}
\hat{R}_L = \sqrt{V(Q_L)C_L}+ \sum\limits_{\ell=0}^{L-1} \sqrt{c_2  M^{-\beta}_\ell\mathcal{C}_\ell},
\end{equation}
and can use this estimate -- in particular $\hat{R}_L^2$ -- to predict if incorporating newer coarse levels will reduce the MLMC cost.}

\section{Numerical results}
\label{section:numresults} 

To demonstrate the utility of redistribution, the following sections compare the performance of multilevel sampling with and without redistribution. Specifically, Section \ref{sec:results1} considers the scaling analysis when solving Darcy's equation as in (\ref{eq:darcy}) on multiple levels of refinement. Subsequently, Section \ref{sec:results2} compares the MLMC performance for a fixed level $\ell=0$ and at least five coarser levels, where redistribution provides the added flexibility of incorporating additional coarser levels into the multilevel hierarchy. 
To achieve scaling on each level in the multilevel hierarchy, we utilize state-of-the-art scalable $H(\div)$ solvers. In particular, the linear system formed to solve the Darcy equation in (\ref{eq:darcy}) is solved with the hybridization AMG approach of ~\cite{lee2017parallel}. Furthermore, \hillary{we apply the PDE-based approach of ~\cite{Osborn17,Osborn17b} to form} realizations of the Gaussian random field $u$; however, to streamline the presentation, the numerical results focus solely on computational time to solve the Darcy equation given an independent realization of the coefficient $k=\exp(u)$. The numerical simulations were completed using tools in the Parelag~\cite{parelag} and ParelagMC~\cite{parelagmc} software, where MFEM~\cite{mfem} was utilized to generate the fine grid finite element discretization, Metis graph partitioner~\cite{metis} was used to form coarse agglomerates, and {\it hypre}~\cite{hypre} to handle parallel linear algebra for both redistribution and forward simulations. We note, all computation was completed on the Quartz cluster at Lawrence Livermore National Laboratory, consisting of 2,688 nodes where each node has two 18-core Intel Xeon E5-2695 CPUs. For the computational results, \hillaryNEW{we use 32 MPI processes per node except for cases with 8 or fewer cores.}

\subsection{Scalability analysis of forward solves}\label{sec:results1}
In this section, scaling analysis is performed on the discrete form of (\ref{eq:darcy}) using the unit cube as the spatial domain. The finest level of the hierarchy is discretized with uniform hexahedral elements such that the local number of elements is approximately \hillaryNEW{$32$K} elements, and, when using $\text{nc}=512$ cores, the global number of elements is about \hillaryNEW{$16$M}. The coarser levels in the multilevel hierarchy are generated by algebraically coarsening by a factor of \hillaryNEW{$8$ elements} (on each level) using the Metis graph partitioner~\cite{metis}. See Table~\ref{tab:ne} for the average number of elements on each level for the $\text{nc}=512$ setting, including the local number of elements. 

When using redistribution, we can {improve the scaling and} incorporate new coarse levels. 
\hillaryNEW{In this work, we compare the standard {\it no redistribution} case with the {\it redistribution} case, both based on the hierarchy of global elements per level provided in Table~\ref{tab:ne}. The standard (no redistribution) case has a maximum of six levels with $\text{nc}=512$ on each level. In the redistribution case, redistribution is utilized to form additional levels in the hierarchy based on the coarsening factor $\beta_c=8$ (see \eqref{eq:core-coarsening-factor}), maintaining no less than 64 elements per core, except when the global number of elements is less than 64}. For these settings, the average number of elements per core and the number of active cores are provided in Table~\ref{tab:ne}. 
\begin{remark}
In this work we do not attempt to identify the optimal number of cores to utilize on each level. Rather, these results are meant to show that redistribution can improve scaling efficiency in MLMC. \hillaryNEW{The value $\beta_c=8$ was selected to match the mesh coarsening factor, thus maintaining a constant local problem size.}
\end{remark}

\begin{table} 
\caption{Average local problem sizes (number of elements, in short NE) for the the standard and redistribution settings, where the number of total cores are provided in parenthesis. }\label{tab:ne}
\begin{center}
\begin{tabular}{p{1cm}p{2.5cm}p{3cm}p{3cm}}
\toprule
     Level $\ell$ & Global NE \newline $\text{nc}=512$ & Local NE (nc) \newline No Redistribution  & Local NE (nc) \newline Redistribution   \\
     \midrule
      0 & 16,777,216 & 32,768 (512) & 32,768 (512) \\
      1 & 2,097,152  & 4,096 (512)  & 4,096 (512)  \\
      2 & 262,144    & 512 (512)    & 512 (512)     \\
      3 & 32,768     & 64  (512)    & 64 (512)     \\
      4 & 4,096      & 8   (512)    & 64 (64)        \\
      5 & 512        & 1   (512)    & 64 (8)      \\
      6 & 64         & -            & 64 (1)       \\
      7 & 8          & -            & 8 (1)       \\
      \bottomrule
\end{tabular}
\end{center}
\end{table}

\hillaryNEW{For the weak scaling analysis, we seek to show that, when changing the global problem size and the number of resources (cores) at the same rate, a constant walltime is maintained. Table~\ref{tab:ne} provides the The local problem sizes on each level for the standard and redistribution case. In this work, we coarsen the global problem size and reduce the number of total cores by a factor of eight, such that the local problem sizes on each level are fixed. We repeat this process, resulting in three problem sizes: {\it Problem 1}, {\it Problem 2}, and {\it Problem 3} with $\text{nc}=8$, $\text{nc}=64$, and $\text{nc}=512$, respectively.}

\begin{figure}[h]
	\centering
	\includegraphics[width=.45\textwidth]{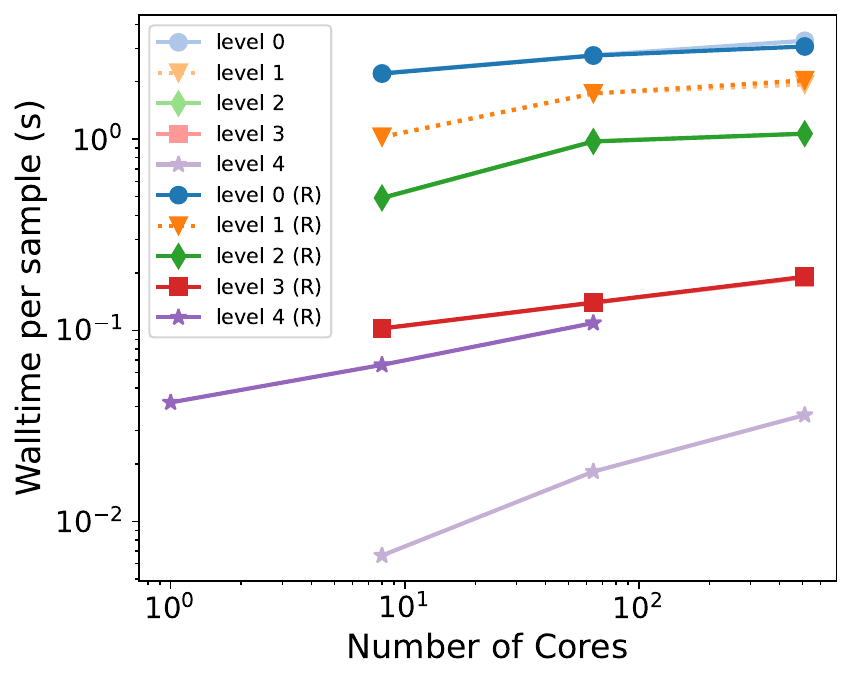}
 	\includegraphics[width=.43\textwidth]{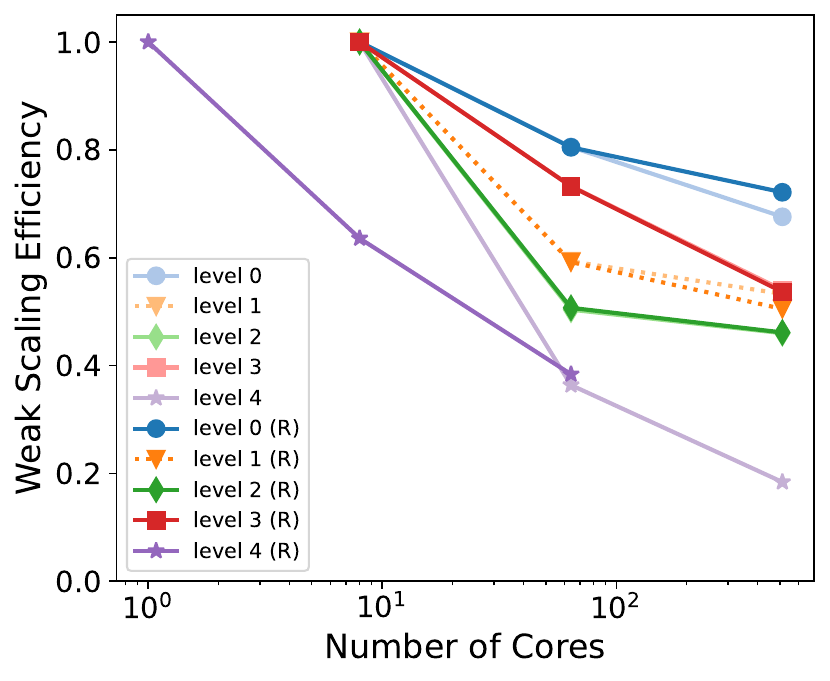}\\
   \hspace*{10pt}(a)\hspace*{170pt}(b)
	\caption{(a) Weak scaling results comparing average walltime of a single Darcy simulations for five levels, where level $\ell=4$ scaling is displayed for simulations with and without redistribution. (b) Weak scaling efficiency for multiple levels as determined by computational time. While, at each level $\ell$, the non-redistribution and redistribution cases share approximately the same global number of elements (see Table~\ref{tab:ne}), the act of redistributing the problem to fewer cores increases the local problem size. The simulation times are provided in Table~\ref{tab:ws-time}.
	 }\label{fig:wstime}
\end{figure}

\hillaryNEW{Table~\ref{tab:ws-time} provides the average walltime to simulate the discretized Darcy's equations for the three problems. The first set of timings, levels $0$-$5$, provide the walltimes for the standard case, while the second set of timings, levels $0$ (R)-$7$ (R), provide those of the redistribution case. Across the three Problems, levels $\ell=0, 1,2,3$ (with and without redistribution) maintain the same local problem size. For levels $\ell=4,5,6,7$, the timings at each level and between the cases are based on different local problem sizes, where the number of cores are provided in parentheses. Figure~\ref{fig:wstime} (a) displays the computational time (walltime averaged over 10 simulations) to simulate the discretized Darcy's equations for the first five levels of each setting. These results allow us to compare the scaling across all levels and between the redistribution cases. As expected, the scaling results for levels $0$-$3$ are approximately the same; however, the levels that utilize redistrbution (level $\ell=4$), show an increase in walltime relative to their non-redistributed counterpart. 
This is due to the changes in local problem size due to redistribution. By performing redistribution, the local problem size increases, and thus the cost increases. What we are more interested in is the slope, which indicates the scaling of the method. Figure~\ref{fig:wstime} (b) provides this as the weak scaling efficiency, which is the computational time is taken with respect to the smallest problem, Problem 1 ($\text{nc}=8$). The key result here is that, with redistribution, level $\ell=4$ shows an improved scaling. On the finest level of Problem 3 ($nc=512$), what was an efficiency of $20\%$ is now an efficiency of $40\%$. }

\begin{table}[h]
\caption{Average wall time (seconds) to solve (\ref{eq:darcy}) using hybridization AMG (averaged over $10$ realizations), for the three problem sizes, shown for the standard core distribution (Level $0$-$5$), and redistribution case denoted (Level 0 (R)-7 (R)). For each timing, the number of utilized cores, nc, is provided in parentheses.}\label{tab:ws-time}
\begin{center}
\begin{tabular}{lllll}
\toprule                                                
Level $\ell$ & Local NE & Problem 1 & Problem 2 & Problem 3 \\ 
&& Walltime (nc) & Walltime (nc)&Walltime (nc)\\ 
\midrule
0 & 32768 & 2.217 (8) & 2.753 (64) & 3.28 (512)\\
1 & 4096 & 1.033 (8) & 1.741 (64) & 1.938 (512)\\
2 & 512 & 0.491 (8) & 0.976 (64) & 1.068 (512)\\
3 & 64 & 0.102 (8) & 0.139 (64) & 0.189 (512)\\
4 & 8 & 0.007 (8) & 0.018 (64) & 0.036 (512)\\
5 & 1 & 0.002 (8) & 0.006 (64) & 0.012 (512)\\
\midrule
0 (R) & 32768 & 2.208 (8) & 2.744 (64) & 3.061 (512)\\
1 (R) & 4096 & 1.028 (8) & 1.737 (64) & 2.035 (512)\\
2 (R) & 512 & 0.494 (8) & 0.974 (64) & 1.071 (512)\\
3 (R) & 64 & 0.102 (8) & 0.14 (64) & 0.191 (512)\\
4 (R) & 64 & 0.042 (1) & 0.066 (8) & 0.109 (64)\\
5 (R) & 8, 64, 64 & 0.002 (1) & 0.013 (1) & 0.031 (8)\\
6 (R) & 8, 64 & - & 0.004 (1) & 0.032 (1)\\
7 (R) & 8 & - & - & 0.005 (1)\\
 \bottomrule   
\end{tabular}
\end{center}
\end{table}

\begin{figure}[h]
	\centering
	\includegraphics[width=.45\textwidth, trim = 0mm 0mm 0mm 0mm]{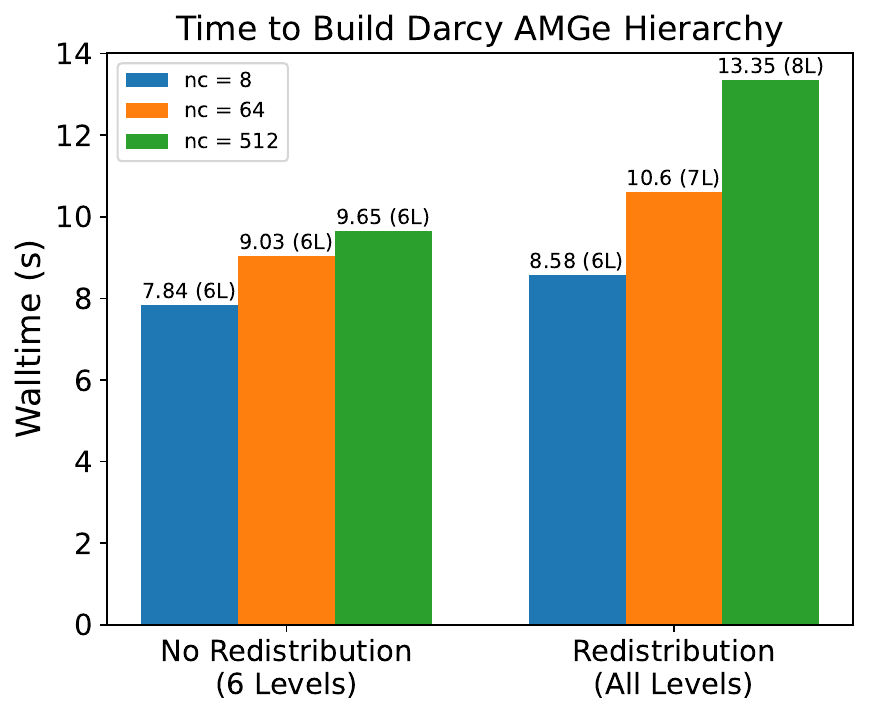}\\
	 \caption{Comparison of the total walltime to build the Darcy AMGe hierarchy for the three Problems, split into the non-redistribution case and the redistribution case. Each bar is labeled with the walltime as well as the number of levels in the hierarchy, e.g., (6L) indicates six levels.}\label{fig:hier}
\end{figure}

\hillaryNEW{In Figure~\ref{fig:hier} we report the total walltime (in seconds) to build the AMGe hierarchy for the Darcy solves. This is done for the three Problems, and for the non-redistribution case and the redistribution case. Incorporating redistribution typically increases the construction time of the AMGe hierarchy (a one-time, offline cost to our MLMC) because more levels are constructed and the data is redistributed to a smaller set of cores on coarser levels.}
\hillaryNEW{When comparing the standard case with the redistribution case, we see no more than a five second increase in the time to build the hierarchy with this new redistribution capability. These timings include building the finest level, which typically takes the most time (with respect to all levels). We do note, the redistribution configuration can have a greater impact on timings depending on the configuration; that is, the cost incurred is related to the number of elements being redistributed. Since the redistribution is not applied until reaching a very small threshold (i.e., a very coarse level), an almost negligible increase in walltime is incurred.}

\subsection{Multilevel Monte Carlo Results}\label{sec:results2}
Using the same hierarchy of levels from the previous section (with the $\text{nc}=512$ case) provided in Table~\ref{tab:ne} (similarly, Problem 3 in Table~\ref{tab:ws-time}), we compare the MLMC performance with and without redistribution \hillaryNEW{when targeting an MSE tolerance of $\varepsilon^2=2.5\mathrm{e}{-5}$}. In particular, we require that the redistribution MLMC approach results in approximately the same statistical estimates as the standard MLMC approach, but with a faster computation time. Recall, the redistribution approach allows for additional coarser levels in the hierarchy. Due to this, it is expected that the total computation time (CPU time) of the MLMC approach with redistribution will be faster than the standard MLMC approach which is limited to six levels.

\begin{figure}
	\centering
	\includegraphics[width=.45\textwidth]{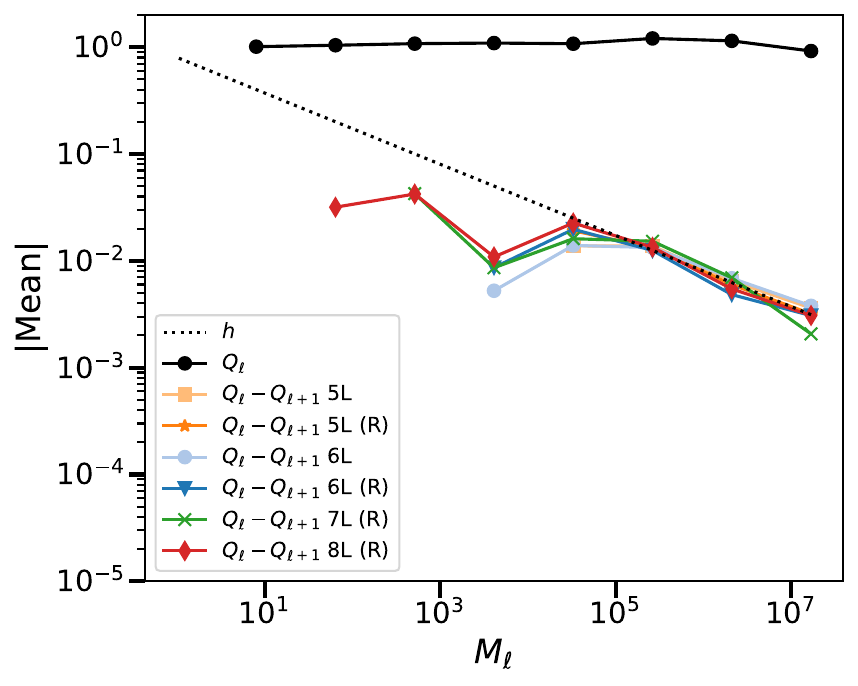}
         \includegraphics[width=.45\textwidth]{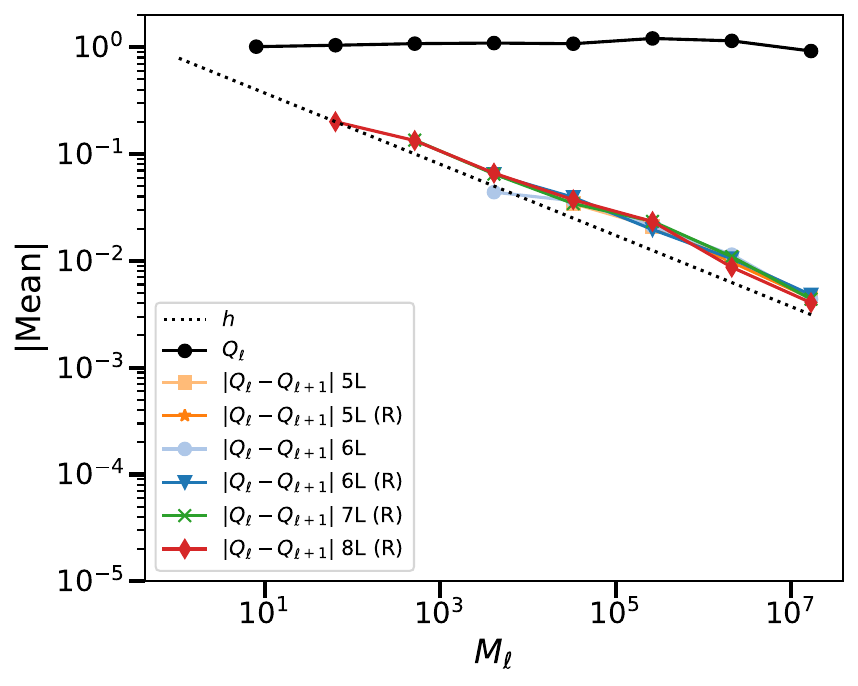}
           \hspace*{35pt}(a)\hspace*{195pt}(b)\\
        \includegraphics[width=.45\textwidth]{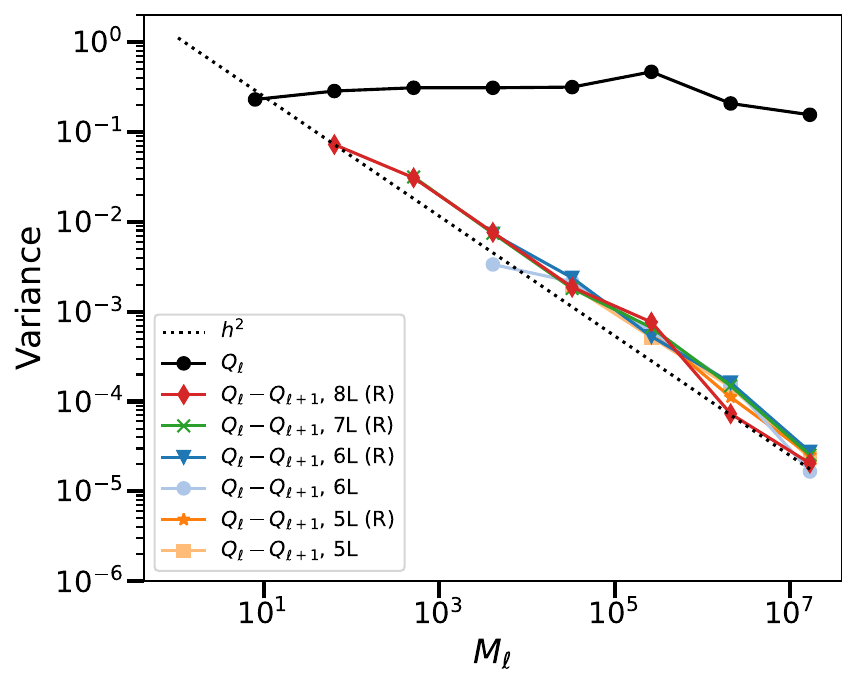}
        \includegraphics[width=.45\textwidth]{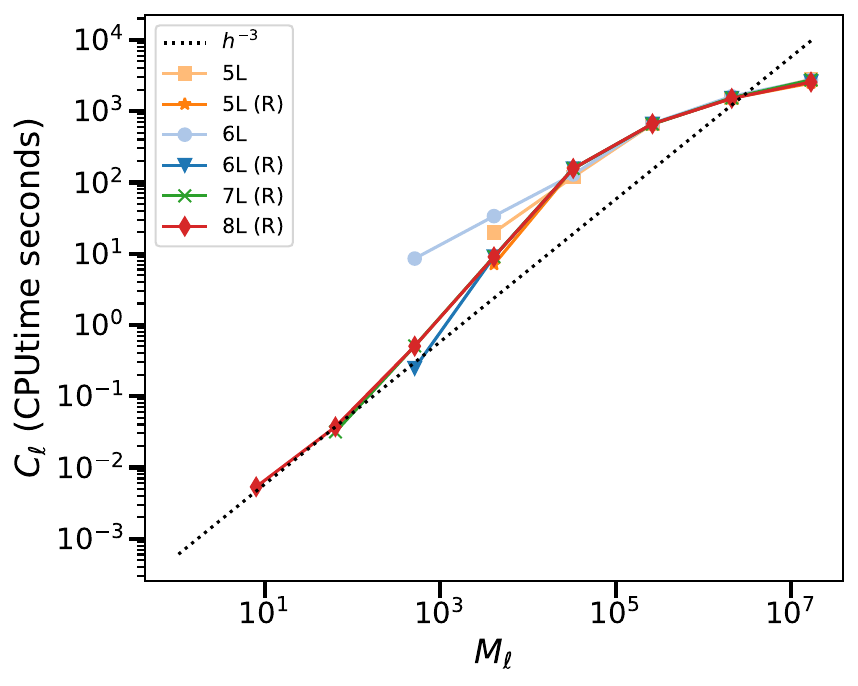}\\
  \hspace*{35pt}(c)\hspace*{195pt}(d)
	 \caption{Comparison of six MLMC settings for the largest problem size ($nc=512$), where the legend indicates the number of levels (e.g. 5L is five levels) and if redistribution was used (``(R)''). (a) Mean approximations of $Q_\ell$ and $Q_\ell-Q_{\ell+1}$ for the various levels in the MLMC hierarchy. (b) Mean approximations of $\vert Q_\ell\vert$ and $\vert Q_\ell-Q_{\ell+1}\vert$ for the various levels in the MLMC hierarchy. (c) Variance approximations of $Q_\ell$ and $Q_\ell-Q_{\ell+1}$ for the various levels in the MLMC hierarchy, and for all settings. (d) Cost on each level (in CPU time). Dashed lines provide comparisons with powers of the mesh size $h$, with the relation $h^{-3}=M_\ell$. }\label{fig:EY}
\end{figure}

The first results considered are the mean and variance estimates of $Q_\ell$ and $Y_\ell=Q_{\ell}-Q_{\ell+1}$. Figure \ref{fig:EY} (a) provides mean estimates, which can be utilized for both the MC mean and MLMC mean estimates (as in (\ref{eq:slmc}) and (\ref{eq:QML})); furthermore, the decay in $Y_\ell$ is useful to assess whether or not the finest level in the hierarchy meets the required discretization error (related to a mean squared error tolerance). \hillaryNEW{Here we compare both settings -- using hierarchies of five to eight levels -- to confirm similar slopes in the $Y_\ell$ estimate, and estimate a value of $\alpha\approx 1/3$ via a log-linear fit between $E[Y_\ell]\approx c_1 M_\ell^{-\alpha}$.} Figure \ref{fig:EY} (b) provides the variance estimates, which are essential to calculating the optimal number of samples to compute on each level in the hierarchy, as defined in (\ref{eq:nl}). In particular, the decay in $\Var(Y_\ell)$ with the refinement of levels indicates that fewer fine samples will be required, thus resulting in a cost reduction compared to single level MC. \hillaryNEW{The notable result here, with respect to the redistribution capability, is the new coarse levels maintain a variance of $Y_\ell$ smaller than that of $Q_\ell$ for $\ell=5$ and $\ell=6$ (note $\ell=7$ is the coarsest level of the eight-level hierarchy, with $Y_7=Q_7$ by definition). While we must still inspect the cost $C_\ell$ and associated number of samples $N_\ell$ on each level to determine if we have cost savings, this result shows promise of MLMC speed up. All settings result in similar variance decays, with $\beta\approx 2/3$}.     

\hillaryNEW{Next we compare the Darcy CPU cost estimates $C_\ell = nc_\ell * (\text{walltime})$ on each level in Figure~\ref{fig:EY} (d). The standard setting, which is limited to six levels, indicates a deterioration in the scaling on the coarsest levels, where the ideal cost scales as $M_\ell = h^{-3}$ (with mesh size $h$). In comparison, the redistribution case improves the coarse level scaling, such that that cost estimates are closer to the $h^{-3}$ line. In this case, we have the approximate relation of $C_\ell \approx c_3 M_\ell^{\gamma}$ with $\gamma=1$. }

\begin{figure}
	\centering
	\includegraphics[width=.45\textwidth]{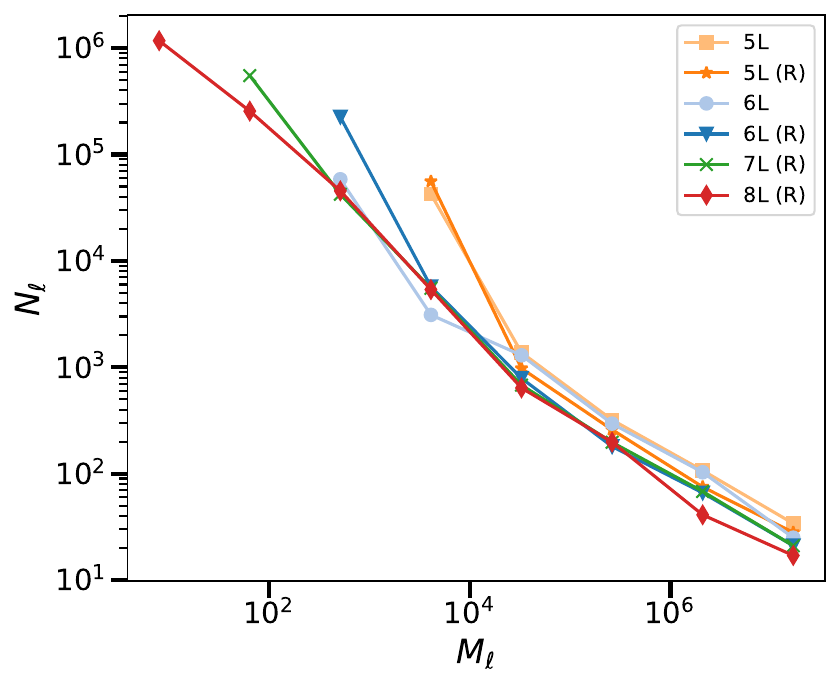}
   \hspace*{15pt}
        \includegraphics[width=.45\textwidth]{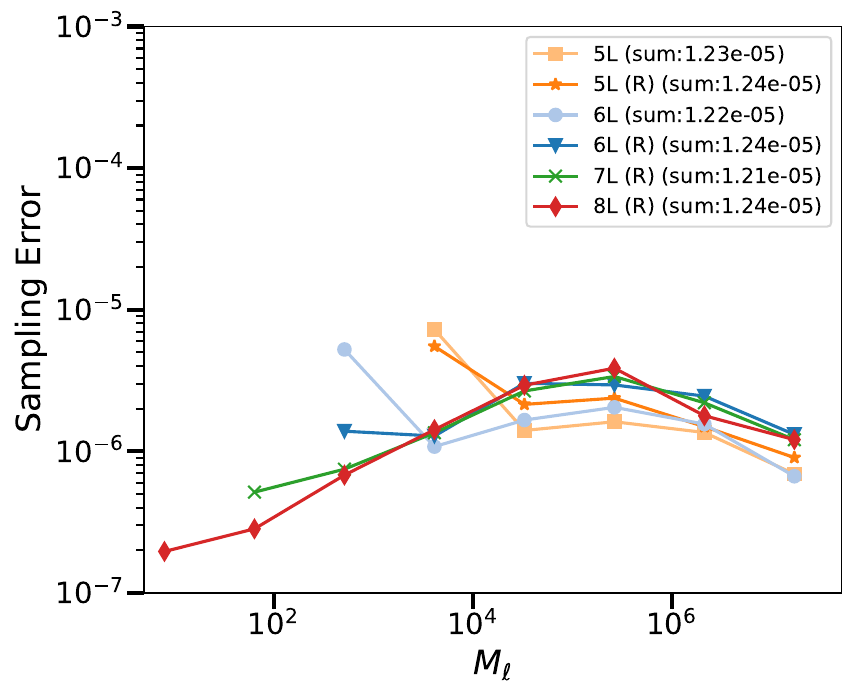}\\
           \hspace*{10pt}(a)\hspace*{170pt}(b)
	 \caption{(a) The number of samples per level (as defined in (\ref{eq:nl})) needed to perform MLMC for the five, six, seven, and eight level cases to attain a sampling error tolerance of $\varepsilon^2/2=1.25\mathrm{e}{-5}$. Note the ``(R)'' indicates redistribution case. (b) The sampling error contributed from each level within the MLMC hierarchy, for each setting. The total sampling error (sum) of these is provided in the legend. }\label{fig:nl}
\end{figure}
Utilizing both the estimated computational time $\mathcal{C}_\ell$ and estimated variance of $Y_\ell$, we estimate the optimal number of samples needed on each level (see (\ref{eq:nl})). 
Figure~\ref{fig:nl} (a) compares the optimal values of $N_\ell$, for each level in the hierarchy, \hillaryNEW{for six cases: five to six level hierarchies for the standard case, and five to eight level hierarchies for the redistribution case.}
As demonstrated, by incorporating more coarse levels in the multilevel hierarchy, fewer simulations \hillaryNEW{from level $\ell=4$ (the coarsest level in the standard five level setting)} are needed. In the redistribution settings, the workload that was solely on $\ell=4$ (or $\ell=5$ in the six level case) is now also distributed amongst new coarser level(s). \hillaryNEW{Figure~\ref{fig:nl} (b) provides the estimated sampling error ($V(Y_\ell)/N_\ell$) on each level, to confirm that the sampling error is bounded by $\varepsilon^2/2=1.25\mathrm{e}{-5}$. This result also shows that, with the addition of coarser levels (as in the redistribution cases), we are able to better balance the error amongst all levels in the hierarchy.}

\begin{figure}
	\centering
	\includegraphics[width=.45\textwidth]{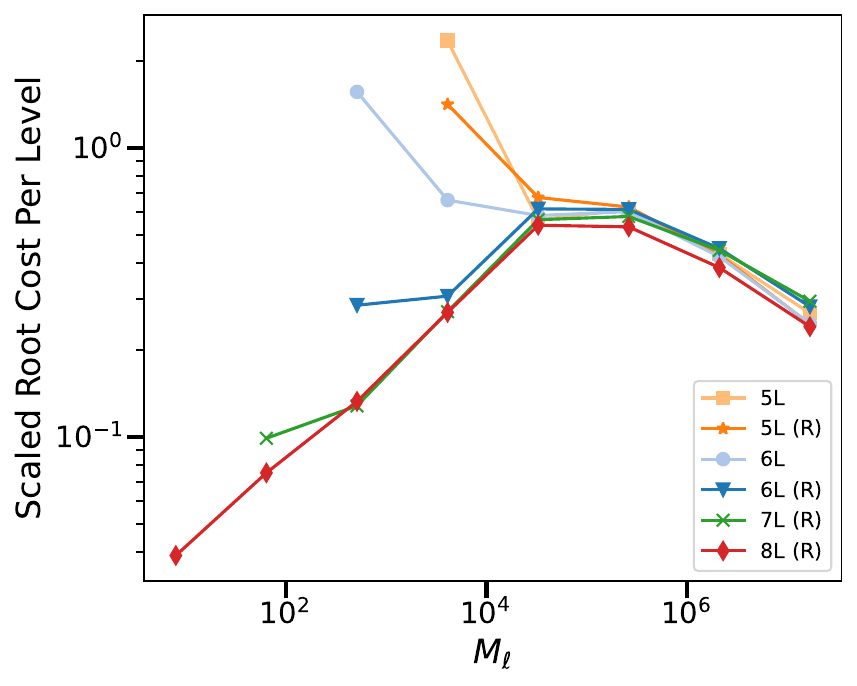}
    \hspace*{15pt}
	\includegraphics[width=.45\textwidth, trim = 0mm 0mm 0mm 0mm]{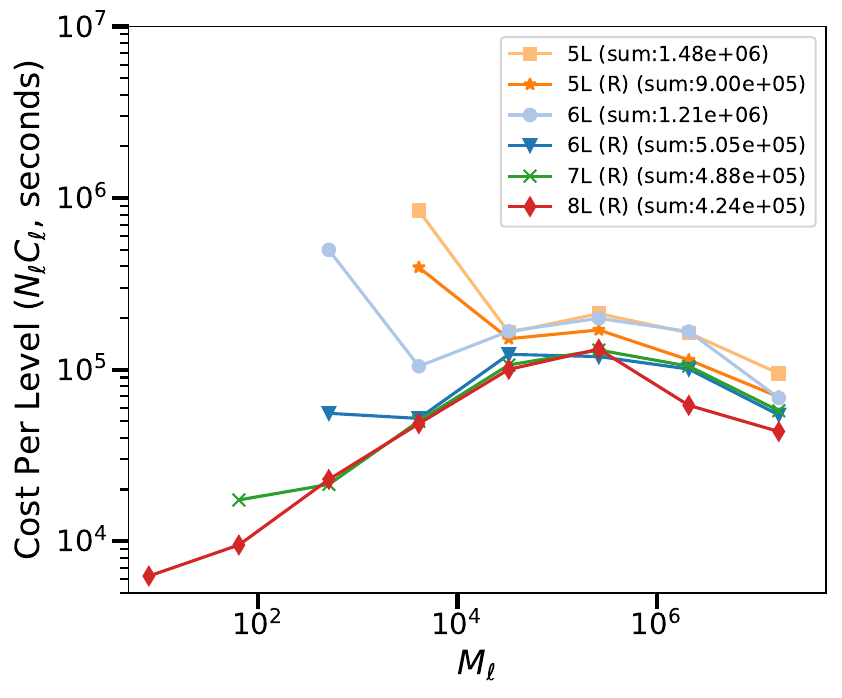}\\
 
	 (a)\hspace*{150pt}(b)
	 \caption{(a) Predicted Scaled Root Cost (see \eqref{eq:scaledrootcost}) on each level for each setting. (b) The total cost to perform MLMC Darcy simulations on each level, for each setting. While results in (a) and (b) have different scales, they do share similar cost distributions amongst the levels. The total MLMC cost \eqref{eq:MLMCcost} is provided in the legend. }\label{fig:mlmccost}
\end{figure}

\hillaryNEW{The cost per level, as well as the total MLMC cost, of the six approaches is provided in Figure~\ref{fig:mlmccost}. Figure~\ref{fig:mlmccost} (a) compares the predicted Scaled Root Cost (see \eqref{eq:scaledrootcost}) on each level for each setting. While this does not provide the true cost on each level, it does indicate how the costs relate. In this case we see the standard MLMC is predicted to have the highest cost incurred on its coarsest level compared to costs on all other levels.} \hillaryNEW{Figure~\ref{fig:mlmccost} (b) shows the terms of total cost sum, $N_\ell \mathcal{C}_\ell$, for each level. This supports the predicted cost relationship between the methods, in particular, it indicates the majority of the MLMC cost comes from the coarsest levels for the standard (non-redistributed) cases. In contrast, the redistribution cases (denoted with (R) in the legend) have a better balanced cost per level; furthermore, the coarsest level costs are reduced by about a factor of $10$.}

\hillaryNEW{When performing redistribution to further coarsen the hierarchy, we see the total cost is reduced. The total cost of each method, calculated via \eqref{eq:MLMCcost}, is provided in the legend of Figure~\ref{fig:mlmccost} (b). As with standard MLMC results, utilizing more coarse levels of discretization in the multilevel hierarchy results in a smaller overall CPU time for the MLMC algorithm. When applying redistribution we are able to obtain at least a $2\times$ speed-up. }

\subsection{MLMC Cost Analysis at Multiple Scales}\label{sec:results3}
\hillaryNEW{Finally, we compare the cost performance of MLMC across the three problem sizes, with and without redistribution. To do this we consider same three problem sizes of Section ~\ref{sec:results1}, Table~\ref{tab:ws-time}; however, now they share the same number of global elements on the coarsest level. That is, without redistribution, Problem 1 has three levels, Problem 2 has four levels, and Problem 3 has five levels, such the coarsest level for all problems has $512$ global elements. With redistribution, Problem 1 has six levels, Problem 2 has seven levels, and Problem 3 has eight levels, such the coarsest level for all problems has $8$ global elements. We set the MLMC mean square error tolerance $\varepsilon^2$ to 4e-4, 1e-4, and 2.5e-5 for Problem 1 ($nc=8$), 2 ($nc=64$), and 3 ($nc=512$) respectively. }

\begin{figure}
	\centering
	\includegraphics[width=.45\textwidth]{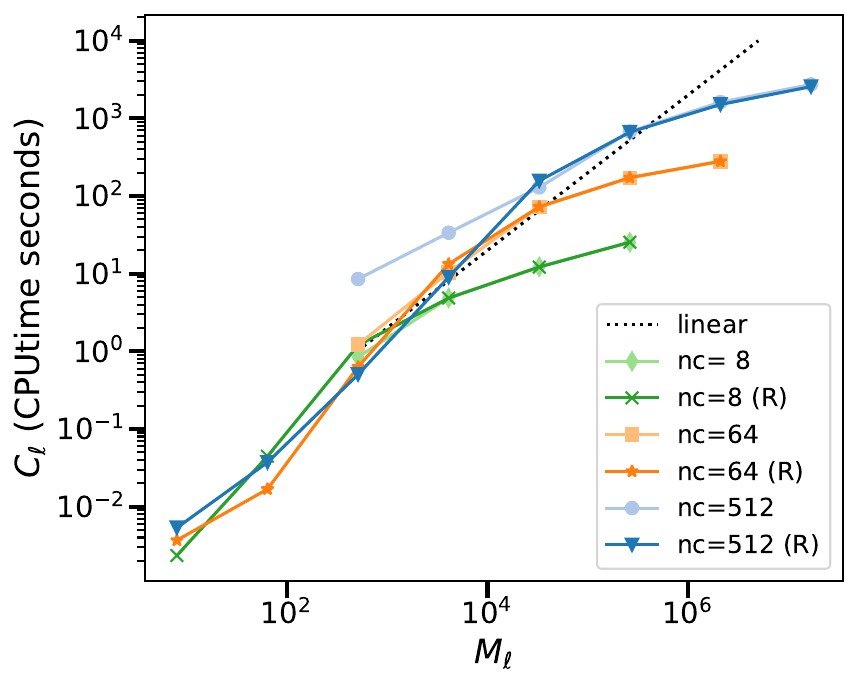}
    \hspace*{15pt}
	\includegraphics[width=.45\textwidth, trim = 0mm 0mm 0mm 0mm]{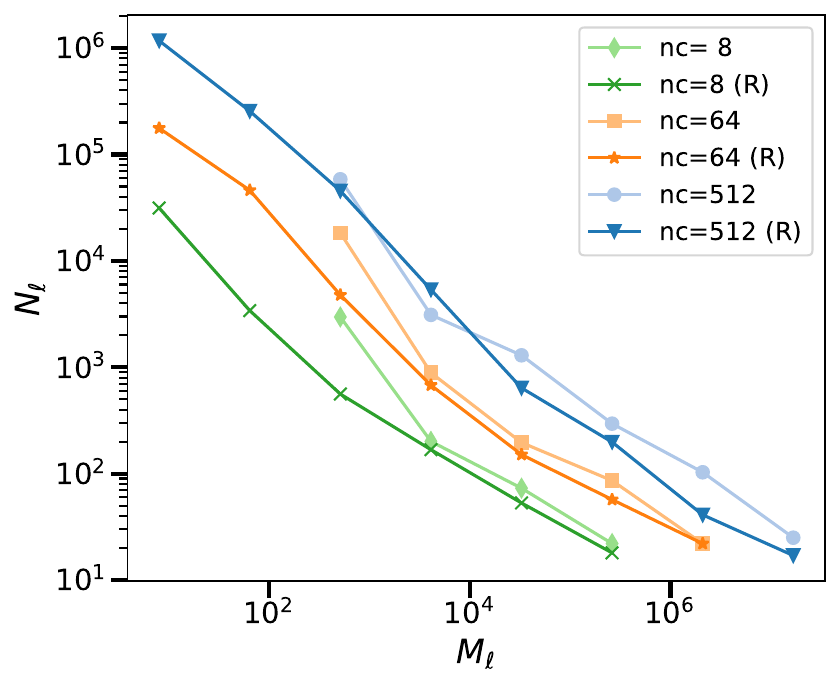}\\

	 (a)\hspace*{150pt}(b)
	 \caption{A comparison of results from each problem size (denoted with $nc=8,64,512$), for the standard case, and the redistribution case (denoted with (R)). (a) The CPU time per sample compared with a linear fit. (b) The optimal number of samples per level needed for MLMC.}\label{fig:all-nl}
\end{figure}

\begin{figure}
	\centering
	\includegraphics[width=.45\textwidth]{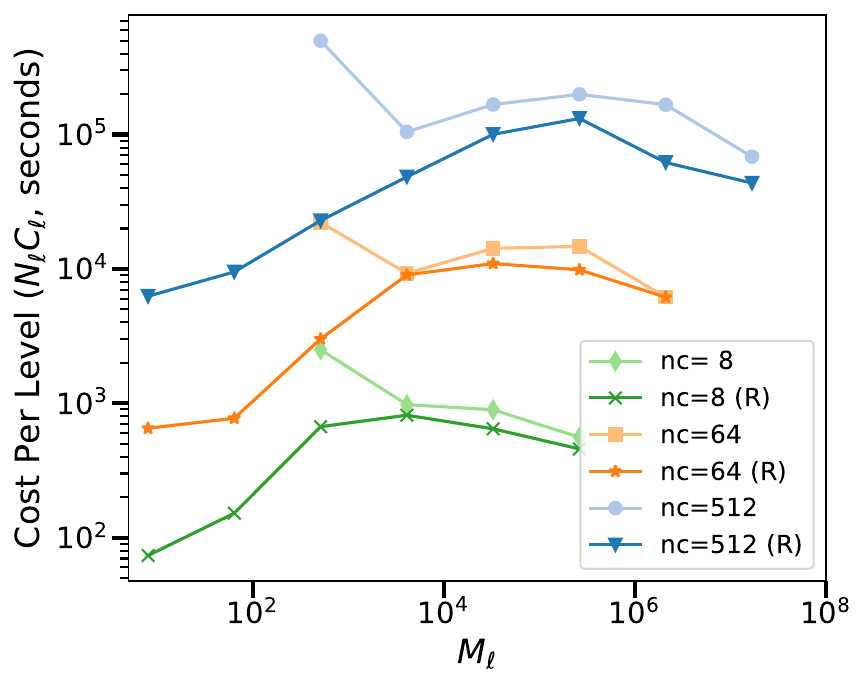}
    \hspace*{15pt}
	\includegraphics[width=.45\textwidth, trim = 0mm 0mm 0mm 0mm]{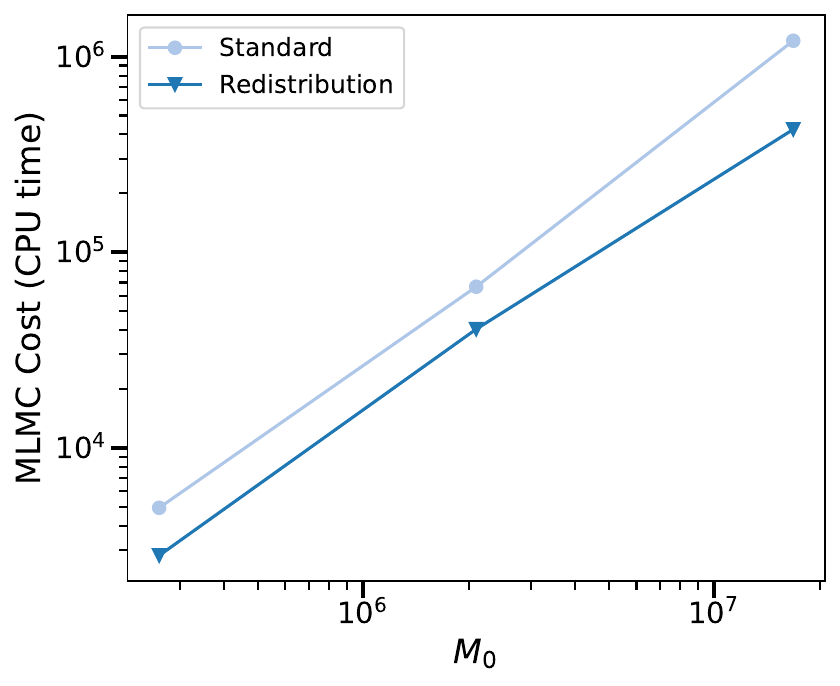}\\
	 (a)\hspace*{150pt}(b) 
	 \caption{A comparison of results from each problem size (denoted with $nc=8,64,512$), for the standard case, and the redistribution case (denoted with (R))(a) CPU cost per level. (b) Total MLMC cost (CPU time) of the standard MLMC compared to the redistributed MLMC for the three problem sizes. Here each problem size is differentiated by the number of elements on the finest level, $M_0$. }\label{fig:all-mlmccost}
\end{figure}

Figure~\ref{fig:all-nl} (a) provides the average CPU time per simulation at each level. The non-redistribution and redistribution cases align fairly closely for changes in mesh size, except for the coarser levels of the $nc=512$ (Problem 3), where the redistribution case clearly performs better. This is expected, as level $\ell=5$ has 1 local element in the non-distribution case and $64$ in the redistribution case. This difference is not as drastic in the $nc=8$ and $nc=64$ cases. The optimal number of samples per level are provided in Figure~\ref{fig:all-nl} (b). Recall, these values are estimated using the cost per level as well as the variance estimates (presented in Section~\ref{sec:results2}). Where we saw, in Figure~\ref{fig:all-nl} (a), Problems 1 and 2 had similar $C_\ell$ estimates between non-redistrbuition and redistribution case, we see the opposite in Figure~\ref{fig:all-nl} (b); in particular, the coarsest levels of the non-redistribution cases require more samples than the associated levels in the redistribution case. Both cases of Problem 3 ($nc=512$) share about the same $N_\ell$ values for levels $\ell=0$ through $\ell=6$. Combining the cost and sample size results provides a clearer picture of the total cost of all approaches, as displayed in Figure~\ref{fig:all-mlmccost} (a).
In each problem size, it is clear that the redistribution is improving the coarse level cost, such that the coarse level is no longer the most computationally expensive. This is because, with the improved scaling, MLMC does a better job balancing the cost across all levels.

The final result, provided in Figure~\ref{fig:all-mlmccost} (b), compares the total CPU time to perform MLMC for the non-redistribution (standard) case and the redistribution case. Here each point is the total cost for given problem size, described by the number of elements on the finest level, $M_0$. The result indicates, that for each problem, utilizing redistribution provides a speed-up. The calculated CPU time speed ups when applying redistribution for Problems 1, 2, and 3 are $1.8\times$, $1.6\times$, and $2.8\times$, respectively.

\section{Conclusions}\label{section: conclusions}
In the present paper, we illustrated an application of parallel redistribution of data at coarse levels for aggregation-based multilevel methods. Such a redistribution is necessary when the number of aggregates becomes smaller than the number of cores used since we must aggregate dofs (and elements, their topological entities, element matrices, etc.) across cores. We summarized the main concepts from \cite{2022ParallelRedistributionReport} in the particular setting of element agglomeration AMG (AMGe) method. A key feature of the studied redistribution is that all operations are formulated in terms of sparse matrix operations (products of matrices and their transposes) performed in parallel. 
Also, we demonstrated that the ability to coarsen beyond the number of cores was also beneficial for the scalability and \hillaryNEW{reduction in CPU time} for MLMC which was our main motivation for the redistribution. The CPU time improvement when using redistribution -- and in particular, for Problem 3 ($nc=512$) relative to Problems 1 and 2 -- suggests that, as these problems move toward extreme-scale, the benefits of utilizing coarse-level redistribution may become more magnified. 
One more feature (not exploited in this paper) is that at coarse levels, since we free some cores, it allows to generate multiple copies of the coarse problems, which \hillaryNEW{will enable an MLMC speed up with regards to walltime; this is in addition to the CPU time performance we demonstrated in this paper}.
This feature and the application of the redistribution in the more complicated multilevel Markov Chain Monte Carlo methods (see, e.g., \cite{Dodwell15} and \cite{2021HierarchicalMLMCMC}) are topics of  %
ongoing research. 

\section*{Acknowledgements}
This work was performed under the auspices of the U.S. Department of Energy by Lawrence Livermore National Laboratory under contract DE-AC52-07NA27344 (LLNL-JRNL-836297). This document was prepared as an account of work sponsored by an agency of the United States government. Neither the United States government nor Lawrence Livermore National Security, LLC, nor any of their employees makes any warranty, expressed or implied, or assumes any legal liability or responsibility for the accuracy, completeness, or usefulness of any information,
apparatus, product, or process disclosed, or represents that its use would not infringe privately owned rights. Reference herein to any specific commercial product, process, or service by trade
name, trademark, manufacturer, or otherwise does not necessarily constitute or imply its
endorsement, recommendation, or favoring by the United States government or Lawrence
Livermore National Security, LLC. The views and opinions of authors expressed herein do not
necessarily state or reflect those of the United States government or Lawrence Livermore
National Security, LLC, and shall not be used for advertising or product endorsement purposes. 

\vspace{0.25cm}
This work of Vassilevski benefited from the Research Training Group (RTG) activities under NSF grant DMS-2136228.%

\bibliographystyle{plainurl}
\bibliography{references.bib}

\end{document}